\documentclass[numbers]{ima-arxiv}

\usepackage{amsfonts}
\usepackage{graphicx}
\makeatletter
\renewcommand*\l@subsubsection[2]{{}{}{}}
\makeatother
\usepackage{subcaption}
\usepackage{epstopdf}
\usepackage{siunitx}  
\usepackage{hyperref}
\usepackage{booktabs}
\usepackage{multirow} 
\usepackage{mathtools} 
\usepackage{amssymb} 
\usepackage{setspace}

\DeclareMathOperator{\erf}{erf}
\DeclareMathOperator{\KL}{KL}
\newcommand{\e}{\mathrm{e}}
\newcommand{\tT}{\mathsf{T}} 
\DeclareMathOperator*{\argmin}{arg\,min}
\DeclareMathOperator{\Cov}{Cov}
\newcommand{\R}{\mathbb{R}}
\newcommand{\N}{\mathbb{N}}
\newcommand{\E}{\mathbb{E}}
\newcommand{\NN}{\mathcal{N}} 
\newcommand{\LL}{\mathcal{L}}
\newcommand{\PP}{\mathcal{P}}

\renewcommand{\d}{\mathrm{d}}

\DeclareMathOperator{\tr}{tr}

\theoremstyle{thmstyletwo}
\newtheorem{theorem}{Theorem}[section]
\newtheorem{proposition}[theorem]{Proposition}
\newtheorem{lemma}[theorem]{Lemma}

\newtheorem{remark}[theorem]{Remark}

\numberwithin{equation}{section}

\newenvironment{romanenumerate}
{%
    \renewcommand{\theenumi}{\roman{enumi}}%
    \begin{enumerate}[(iii)]%
}
{%
    \end{enumerate}%
}

\makeatletter
\newcommand{\beginsupplement}{%
  \clearpage

  \setcounter{section}{0}%
  \setcounter{figure}{0}%
  \setcounter{table}{0}%
  \setcounter{equation}{0}%

  \renewcommand{\thesection}{S\arabic{section}}%
  \renewcommand{\thefigure}{S\arabic{figure}}%
  \renewcommand{\thetable}{S\arabic{table}}%
  \renewcommand{\theequation}{S\arabic{equation}}%

  \renewcommand{\theHsection}{supp.section.\arabic{section}}%
  \renewcommand{\theHfigure}{supp.figure.\arabic{figure}}%
  \renewcommand{\theHtable}{supp.table.\arabic{table}}%
  \renewcommand{\theHequation}{supp.equation.\arabic{equation}}%
}
\makeatother

\hypersetup{
  pdftitle={A Jump-Diffusion Framework for Irregular Time Series Generation},
  pdfauthor={Onno Pfohl, Jannis Chemseddine, Paul Hagemann, Gabriele Steidl, Christian Wald, Tim Jahn},
  pdfkeywords={generative modeling, generator matching, jump-diffusion processes, time series}
}

\begin{document}

\firstpage{1}
\title[A Jump-Diffusion Framework for Irregular Time Series Generation]{A Jump-Diffusion Framework for Irregular Time Series Generation}

\author{Onno Pfohl* \ORCID{0009-0006-7546-3028}
\address{\orgdiv{Department of Mathematics}, \orgname{Humboldt-Universität zu Berlin}, \orgaddress{\postcode{12489}, \state{Berlin}, \country{Germany}};
and
\orgdiv{Department of Mathematics}, \orgname{Technische Universität Berlin}, \orgaddress{\postcode{10623}, \state{Berlin}, \country{Germany}}}}

\author{Jannis Chemseddine \ORCID{0009-0009-7460-6367} 
\address{\orgdiv{Department of Mathematics}, \orgname{Technische Universität Berlin}, \orgaddress{\postcode{10623}, \state{Berlin}, \country{Germany}}}}

\author{Paul Hagemann \ORCID{0009-0006-9679-5654}
\address{\orgname{Bundesanstalt f\"ur Materialforschung und -pr\"ufung}, \orgaddress{\postcode{12205}, \state{Berlin}, \country{Germany}}}}

\author{Gabriele Steidl \ORCID{0000-0003-4638-9245}
\address{\orgdiv{Department of Mathematics}, \orgname{Technische Universität Berlin}, \orgaddress{\postcode{10623}, \state{Berlin}, \country{Germany}}}}

\author{Christian Wald \ORCID{0000-0003-2373-9492}
\address{\orgdiv{Institut Camille Jordan, UMR 5208},
\orgname{INSA Lyon, Centrale Lyon, Université Lyon 1, Université Jean Monnet, CNRS},
\orgaddress{
\postcode{69621},
Villeurbanne,
\country{France}
}}}

\author{Tim Jahn \ORCID{0000-0002-7633-8265}
\address{\orgdiv{Department of Mathematics}, \orgname{Technische Universität Berlin}, \orgaddress{\postcode{10623}, \state{Berlin}, \country{Germany}}}}
\authormark{O. Pfohl et al.}

\corresp[*]{Corresponding author: \href{email:onno.pfohl@hu-berlin.de}{onno.pfohl@hu-berlin.de}}

\abstract{We propose a framework for generative modeling of continuous-time processes from irregularly and asynchronously recorded data. It is based on the matching of generators and accommodates discontinuous trajectories. Analytical formulas for diffusion and jump bridges yield a family of reference generators that a neural network is trained to match. The key ingredient is that, for our constructed jump bridge, a parametrization of the jump kernel densities by scaled Gaussians admits closed-form expressions for the Kullback-Leibler divergence, allowing simulation-free training.
}


\maketitle

\section{Introduction}
Generative modeling has shown impressive results in imaging and natural sciences \cite{corso2023diffdock, rombach2022high, pmlr-v37-sohl-dickstein15, song2021scorebased}. However, the generation of time series, which arise in the context of modeling financial, meteorological and patient vital data has received less attention \cite{buehler2020generating, jia2019neural, kerrigan2024functional, kidger2021neural, zhang2024trajectory}. 
In particular, in the study of financial markets and price prediction \cite{GULMEZ2023120346,wu2023graph}, the (conditional) generation of time series data is an important cornerstone. For instance, in deep hedging \cite{buehler2019deep}, simulation of market data is crucial to learn appropriate hedging strategies.

The usual approaches to model the time series distribution involve adversarial training \cite{kidger2021neural},  backpropagation through the SDE solver \cite{jia2019neural}, or fitting an object as high dimensional as the number of discrete time points \cite{kerrigan2024functional}.  Another difficulty  many methods face is that they rely on  fixed time discretizations, and therefore are not feasible e.g. for many medical datasets, where points are measured irregularly. Although it is sometimes possible to fill those missing values as in \cite{yoon2019time}, we are interested in exploring a method which is not reliant on a fixed time grid and can directly accommodate irregular and varying 
 observation times. Throughout, observations are linked across time into individual trajectories, unlike in trajectory inference \cite{MMFM,MMSFM, accelerationmatching}, where only unpaired samples of each marginal are available.

\emph{Simulation-free} approaches, which do not have to solve an ODE/SDE during training, 
 such as flow matching \cite{albergo2025stochastic,lipman2023flow,liu2023flow, Ws2025} or score-based diffusion \cite{song2021scorebased}, including recent extensions to infinite-dimensional spaces \cite{kerrigan2024functional,hagemann2025multilevel}, and static Schr\"odinger bridge estimation \cite{pooladian2025plug} have significantly improved the performance of generative models compared to adversarial networks \cite{gan} or simulation-based methods \cite{chen2018, tzen2019theoretical}. Recently, such simulation-free methods have been applied to time series generation employing a high-dimensional latent distribution in \cite{bilos2023,kerrigan2024functional}. In contrast, the authors of \cite{zhang2024trajectory} proposed to view time series generation as essentially generating one-dimensional curves via interpolating between the discretization points. This interpolation is done by learning an SDE between the discrete time points and taking the realizations at previous time points as conditions. A very similar idea has been proposed by \cite{chen2024probabilistic} in the language of stochastic interpolants and applied to more challenging high-dimensional datasets, and a similar algorithm is proposed in \cite{elgazzar2025probabilisticforecastingautoregressiveflow}. 
 During inference, the autoregressive generation is akin to doing the next word prediction in GPT models \cite{NIPS2017_3f5ee243}, highlighting the conceptual similarity between 
sequential text generation and sequential time series synthesis, while differing crucially in that the latter must handle 
irregular temporal spacing. Another recent approach \cite{splineflow} replaces linear conditional paths by B-spline interpolants. In contrast to the autoregressive approaches above, the learned dynamics are Markovian.

What these approaches have in common is that they generate continuous paths. Processes in economics, however, are often fundamentally discontinuous as  buying or selling a stock may induce a jump.
Jump processes have been already employed in \cite{campbell2023transdimensional} for the forward process inside a diffusion framework \cite{song2021scorebased}, in \cite{benton} for parameterizing and learning generators on discrete state spaces via Markov chains or investigating jump Schrödinger bridges via h-transforms, see \cite{zlotchevski2025schrodingerbridgeproblemjump}. For static and structured problems, e.g. image and protein generation, learning the generators of jump and diffusion processes was considered already in \cite{holderrieth2025generator} under the name of generator matching.

In this paper, we combine generator matching
with the dynamics of time series generation 
which poses both theoretical and numerical challenges.
Since the approach of \cite{zhang2024trajectory} has a singularity when reaching the final time point, which is detrimental in either training or inference, we modify it by allowing small amounts of Gaussian noise at known discretization points, stabilizing the diffusion bridges. As the main contribution, we derive jump bridges for interpolation, and consider ultimately a combination of both. A common method to approximate jump generators is via binning the state space, which creates numerical challenges. In contrast, in order to reduce computational burden, we parameterize the jumps by Gaussian-like functions. This allows us to derive an analytical expression for the  Kullback-Leibler divergence used in the corresponding loss function, addressing a well-known theoretical bottleneck for continuous-time generative models \cite{conforti2025kl}.
Since we model non-Markovian time series, we need to include memory as conditions in our next point predictors.
We prove that our learned time series approximates the true joint distribution and does not only capture marginals.
We evaluate our method on different data sets and compare it with recent time series generation methods.

\textbf{Contributions.} 
Our work provides a new perspective for modeling and generating autoregressive time series.
We construct a flexible generative model for irregularly sampled time series that allows observation times to vary across trajectories.
Our main contributions are the following:
\begin{itemize}
\item In Theorem \ref{th1}, we introduce jump bridges, and modify the diffusion bridges considered by \cite{zhang2024trajectory} for numerical stabilization on stabilized paths. With this approach we avoid singularities for both jump and diffusion bridges, see Figure \ref{fig:prop2} for an illustration. This is key for making training stable and it increases the method's flexibility to account for discontinuous trajectories. 
\item 
We propose to parametrize the density of the jump kernel by scaled Gaussians and to learn the corresponding parameters. This allows us to derive closed-form expressions for the Kullback-Leibler loss function, as presented in Theorem \ref{prop2}. These expressions can be computed efficiently via Theorem \ref{prop1}, avoiding the expensive binning method used in \cite{holderrieth2025generator} and making training scalable.
\item We establish a stability result, Theorem \ref{rem:joint_post}, that controls the approximation error for the joint marginal distribution of the underlying process at any finite collection of times.
\end{itemize}

Based on the theoretical framework, we design a practical algorithm and validate it through numerical experiments in various dimensions. The emphasis of this paper lies on the rigorous derivation of the underlying mathematical objects, with the experiments demonstrating the algorithm they induce. All proofs are presented in Appendices~\ref{app:proofs} and \ref{app:joint_post}. Further numerical results, implementation details and some additional definitions used throughout the paper are given in the Supplementary Material.

\section{Problem Setting}
\label{sec:problem_setting}
\begin{figure}[H]
\centering
\begin{minipage}[t]{0.5\linewidth}
\vspace{0pt}
    \includegraphics[clip, trim={0.2cm 0.28cm 0.2cm 0.3cm},width=\textwidth]{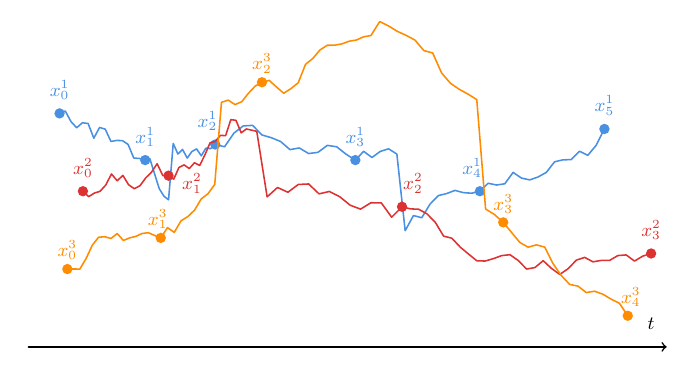}
    \end{minipage}
    \hfill
    \begin{minipage}[t]{0.45\linewidth}
    \vspace{0pt}
    \caption{Time series with irregularly sampled times. Curves: continuous-time càdlàg paths. Big dots: irregularly sampled time-stamps. For readability we write $x^i_j$ instead of $x^i_{t^i_j}$.}\label{fig:irregular_dataset}
    \end{minipage}
\end{figure}
\begin{figure}[H]
\centering
\begin{minipage}[t]{0.5\linewidth}
\vspace{0pt}
  \includegraphics[clip, trim={0.2cm 0.4cm 0.2cm 1cm},width=\textwidth]{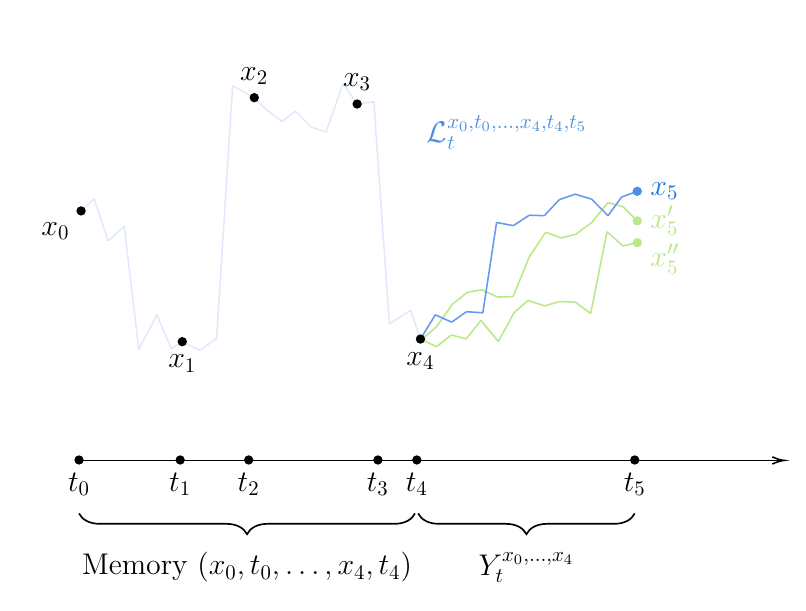}
    \end{minipage}
    \hfill
    \begin{minipage}[t]{0.49\linewidth}
    \vspace{0pt}
    \caption{Process $Y_t^{x_0,\ldots,x_4}$ parametrized by its generator $\mathcal{L}_t^{x_0,t_0,\ldots,x_4,t_4,t_5}$ for $t\in[t_4,t_5]$. The trajectories in $[t_4,t_5]$ are drawn from $Y_t^{x_0,\ldots,x_4}$ and thus $x_5,x'_5,x^{''}_5\sim P_{X_{t_{5}}|X_{t_4}=x_{4},\ldots,X_{t_0}=x_{0}}$. In this illustration out of all possible paths from $Y_t^{x_0,\ldots,x_4}$ we sampled the blue one, and moving forward, would add $x_5,t_5$ to the memory to compute the next process $Y^{x_0,\ldots,x_5}_t$ and predict $x_6$.}\label{fig:conditional_generator}
    \end{minipage}
\end{figure}

We consider a dataset of $N$ time series $\mathcal{X}:=\{x^{1},...,x^{N}\}$, where each time series is defined as $x^i \coloneqq (x^i_{t^i_j})_{j=0}^{n_i} \in \mathbb R^d$. These are sampled at irregular timestamps $0 = t^i_0 < t^i_1 < \ldots < t^i_{n_i} = T$. We assume these sequences are discrete observations of a continuous-time random process $X=(X_t)_{t\in[0,T]}$ taking values in the space of right-continuous functions with left limits (càdlàg). An example of such a dataset can be seen in Figure \ref{fig:irregular_dataset}.
Our objective is to learn the finite-dimensional marginal distributions of $X$, i.e. $P_{X_{t_0},...,X_{t_n}}$, for all finite sequences $0=t_0<t_1<...<t_n= T$, $n\in \N$.  The distribution of $X$ is uniquely determined by all these finite-dimensional marginals \cite[Section 13]{billingsley2013convergence}. The joint distribution $P_{X_{t_0},...,X_{t_n}}$ can be factorized using the disintegration theorem (see Supplementary Material, Section \ref{sec:supp:defs}):
\begin{align}\label{joint_disintegrate}
    P_{X_{t_0},\ldots,X_{t_n}}
    =P_{X_{t_n}|X_{t_{n-1}}=x_{n-1},\ldots ,X_{t_0}=x_0}\times_{x_{n-1},\ldots,x_0}\cdots \times_{x_1,x_0}P_{X_{t_1}|X_{t_0}=x_{0}}\times_{x_0}P_{X_{t_0}}.
\end{align}
This factorization implies an inductive sampling scheme. We first sample an initial value $x_0 \sim P_{X_{t_0}}$. Subsequently, given a history of values $x_0, \ldots, x_j$, the next value $x_{j+1}$ is drawn from the conditional distribution $P_{X_{t_{j+1}}|X_{t_{j}},\ldots,X_{t_0}}$.
In order to sample from the latter distribution, we will learn for 
$x_0,\ldots,x_j\in \R^{d}$ a Markov process $(Y_t^{x_0,\ldots,x_j})_{ t\in[t_j,t_{j+1}]}$  starting in $Y_{t_j}^{x_0,\ldots,x_j}=x_j$ such that approximately
\begin{align} \label{posterior}
    Y_{t_{j+1}}^{x_0,\ldots,x_j} \sim P_{X_{t_{j+1}}|X_{t_j}=x_{j},\ldots,X_{t_0}=x_{0}},
\end{align}
see Figure \ref{fig:conditional_generator}.
Markov processes are convenient because they are analytically tractable. Note that each Markov process interpolates locally between times $t_j$ and $t_{j+1}$ and depends on the previous observations $(x_0,t_0,\ldots, x_j,t_j)$, which we call the memory.

In the following section we explain how to set up and learn the interpolating Markov process for which we employ ideas from generator matching.

\section{Learning the Interpolating Markov Process}
\label{sec:interpolating_mps}
We now detail our approach for approximating the conditional distributions $P_{X_{t_{j+1}}|X_{t_j}=x_{j},\ldots,X_{t_0}=x_{0}}$ by learning interpolating Markov processes. We first characterize these processes via their infinitesimal generators, which serve as the primary objects for neural parametrization. Our main theoretical contribution in this section is the derivation of explicit generators for diffusion and jump bridges, which allow for efficient parameter estimation without relying on nested numerical simulations. In particular, for the jump process generator, we employ a specific parameterization that enables us to compute the KL divergence describing the loss in closed form, significantly reducing computational overhead.

\subsection{Generators of Markov Processes} 
We first introduce Markov processes, generators and the corresponding Kolmogorov forward equation. Let $\mathcal{B}(\R^d)$ denote the Borel sets on $\R^d$ and 
$\mathcal{B}_b(\R^d)$ the space of bounded, measurable functions on $\R^d$. Further, let
$\mathcal M_+(\R^d)$ denote the space of finite non-negative measures on $\R^d$, $\mathcal P(\R^d)$ the probability measures and $\mathcal P_2(\R^d)$ the probability measures with finite second moments.

We consider a  Markov process $(Y_t)_{t\in [0,1]}$ in $\R^d$. The process 
can be  characterized by its generator
$(\mathcal L_t)_t$, 
mapping from a subset of $\mathcal{B}_b(\R^d)$ into the  measurable functions by
\begin{align}
    \mathcal L_t f(x)\coloneqq \lim_{h\downarrow 0}\frac{\E[f(Y_{t+h})-f(Y_t)|Y_t=x]}{h}.
\end{align}
It can be rewritten in the form \cite{ von1965fast}
\begin{align}\label{eq:generator}
    \mathcal{L}_t f(x) = \underbrace{ \langle\nabla_x f(x),u_t(x)\rangle}_{\rm drift} + \underbrace{\frac{1}{2}\langle \Sigma_t^2(x), \nabla_x^2 f(x)\rangle }_{\rm diffusion} + \underbrace{\int(f(y)-f(x)) \, Q_t({\rm d}y,x)}_{\rm jump},
\end{align}
where $u_t:\R^d\to\R^d$ is the drift, $\Sigma_t:\R^d\to\R^{d\times d}$ is the diffusion coefficient and $Q_t:\mathcal B(\R^d)\times\R^d\to\R$ is the rate (jump) kernel (i.e., $x\mapsto Q_t(B,x)$ is measurable and $B\mapsto Q_t(B,x)$ is a finite non-negative measure). Assume that $Y_t$ has a law with density $p_t$. The evolution is described by the \emph{Kolmogorov forward equation} (KFE)
$$
    \partial_t \langle p_t,f \rangle = \langle p_t, \mathcal L_t f \rangle.
$$
We first restrict our attention to $\Sigma_t^2 = \eta ^2 I_d$ with an appropriately given $\eta > 0$ and we assume that $Q_t(\cdot,x)   \in \mathcal M_+(\R^d)$ admits a Lebesgue density $q_t(y,x) \d y$. Then, taking the adjoint of the KFE provides the following PDE 
\begin{equation}
    \partial_t p_t = - \nabla \cdot (p_t u_t) + \frac12 \eta^2 \Delta p_t + \int p_t(y) q_t(x,y) - p_t(x) q_t(y,x) \d y.
\end{equation}

Given a generator $\LL_t$, which we parametrize in the following by the triple $(u_t, \eta, q_t)$, we can approximate a corresponding Markov process by Algorithm \ref{alg:markov_from_generator}. We decompose the density $q_t(\cdot,x) = \lambda_t(x) J_t(\cdot,x)$ into the jump intensity $\lambda_t(x):= \int q_t(y,x){\rm d}y\ge 0$ and a probability density $J_t(\cdot,x)$.
\begin{algorithm}[H]
\caption{Approximating a Markov process}
\label{alg:markov_from_generator}
\begin{algorithmic}[1]
    \Require $u_t, \eta, q_t(\cdot,x) = \lambda_t(x) J_t(\cdot,x)$
    \State \textbf{Initialization:} $x_0 \sim Y_0$, step size $h = 1/n$
    
    \For{$i = 0$ \textbf{to} $n-1$} 
        \State Set $t = i/n$
        \State Sample $x^{J}_{t+h} \sim J_t(\cdot,x_t)$
        \State Sample $m \sim \text{Bernoulli}(h\lambda_t(x_t))$ and $\epsilon_t \sim \mathcal{N}(0,I_d)$
        \State $x^{D}_{t+h} = x_{t} + h u_t(x_t) + \sqrt{h} \eta  \epsilon_t$
        \State $x_{t+h} = m x^{J}_{t+h} + (1-m)x^{D}_{t+h}$
    \EndFor
    
    \Return  $(x_0,x_{\frac{1}{n}},\ldots, x_1)$ approximately distributed as $(Y_0,Y_{\frac{1}{n}},\ldots,Y_1)$
\end{algorithmic}
\end{algorithm}

\begin{remark}
    Note that for $u_t=0$, $\eta=0 $ we have $x^{D}_{t+h}= x_t$ and thus obtain a pure jump process. On the other hand, for $\lambda_t=0$ we have that $m=0$ and thus $x_{t+h}=x^{D}_{t+h}$ is a pure diffusion process.
\end{remark} 

\subsection{Conditional Generators for Interpolation}
We saw that in order to sample from $P_{X_{t_{j+1}}|X_{t_j}=x_{j},\ldots,X_{t_0}=x_{0}}$ in \eqref{posterior}, it suffices to learn the generator of the process $Y_{t_{j+1}}^{x_0,\ldots,x_j}$ and then use Algorithm \ref{alg:markov_from_generator}. To this end, we aggregate simple conditional Markov processes. In this subsection, we derive suitable candidates, whose drift-diffusion and jump generators can be computed analytically.

For notational simplicity we describe the procedure for $P_{X_{t_1}|X_{t_0} =x_0}$, where $x_0 \in \R^d$ is arbitrarily fixed, and assume that $t_0 = 0$ and $t_1 = 1$. The extension to arbitrary $t_0$ and $t_1$ is presented in Appendix \ref{app:proofs}.

We start with a family $(P_t)_{t\in[0,1
]}$ of Markov kernels  $P_t(\cdot, \cdot, \cdot): \mathcal B(\R^d)\times \R^d\times\R^d\to [0,1]$ with
\begin{align}\label{diracs}
    P_{0}(\cdot,x_{0},x_{1}) = \delta_{x_{0}}\quad \text{ and}\quad P_1(\cdot,x_0,x_1)= \delta_{x_1}.
\end{align}
We consider 
$$
\alpha_t^{x_0} 
\coloneqq 
P_t(\d x, x_{0},x_{1})\times_{x_1}P_{X_1|X_0=x_0}(\d x_1) \in \mathcal P_2(\R^d \times \R^d)
\quad \text{and} \quad
P^{x_0}_t\coloneqq \pi^x_\sharp( \alpha_t^{x_0}).
$$
Let $\LL_t^{x_0,x_1}$ be generators of processes $Y^{x_0,x_1}_t$ with marginals $Y^{x_0,x_1}_t\sim P_t(\cdot, x_0,x_1)$. By \eqref{diracs}, $Y^{x_0,x_1}_t$ interpolates between $x_0$ and $x_1$.
Then by \cite[Proposition 1]{holderrieth2025generator} $\mathcal L^{x_0}_t$ defined by 
\begin{equation}\label{int}
    \mathcal L^{x_0}_t(f)(x)
    \coloneqq \int \LL_t^{x_0,x_1}(f)(x)\d (\alpha_{t}^{x_0})_x (x_1)
\end{equation}
generates a  Markov process $Y_t^{x_0}$ with
$Y_t^{x_0} \sim P_t^{x_0}$. By construction $Y_t^{x_0}$ interpolates between $P_0^{x_0}=\delta_{x_0}$ and $P_1^{x_0}=P_{X_1|X_0=x_0}$.
 
\begin{remark}\label{prop:post_x_0}
    Instead of \eqref{diracs}, we will more generally start with Markov kernels $P_t(\cdot,x_0,x_1)$ with 
    $$P_0(\cdot,x_0,x_1) =\NN(x_0,\rho^2 I_d) 
    \quad \text{and} \quad 
    P_1(\cdot,x_0,x_1) =\NN(x_1,\rho^2 I_d), \quad \rho >0$$
    Let $\LL_t^{x_0}$ be the generator from \eqref{int} and let $Y_t^{x_0}$ be an associated Markov process. Then we have $Y_0^{x_0}\sim\NN(x_0,\rho^2 I_d)$ and $Y_1^{x_0}\sim P_{X_1|X_0=x_0}*\NN(0,\rho^2 I_d)$, see Appendix \ref{app:proofs}. 
\end{remark}

Now we construct appropriate $P_t(\cdot, x_0,x_1)$ and $\LL_t^{x_0,x_1}$.

\begin{theorem}\label{th1}
    For arbitrary fixed $x_0,x_1 \in \R^d$,
    set $m_t \coloneqq (1-t)x_0 + t x_1$ and
    $\tau_t\coloneqq \eta^2 t(1-t)+\rho^2$ with
    $\eta > 0$. Let $P_t = P_t(\cdot,x_0,x_1) = \mathcal N(m_t,\tau_t I_d)$ with density $p_t$. Then each of the following generators 
    $\mathcal L_t^{{\rm Diff}} = \mathcal L_t^{{\rm Diff},x_0,x_1}$
    and
    $\mathcal L_t^{\rm{Jump}} = \mathcal L_t^{{\rm Jump},x_0,x_1}$
    fulfills the KFE and thus gives rise to a Markov process $Y_t = Y_t^{x_0,x_1}$ with marginals $Y_t \sim P_t$:
    \begin{romanenumerate}
        \item $\mathcal{L}^{{\rm Diff}}_t f (x) =  \nabla f(x)^\tT u_t(x)  + \frac{\eta^2}{2} \Delta f (x)$, where 
        $u_t^{x_0,x_1} =  u_t \coloneqq  x_1-x_0 - (x -m_t) \frac{\eta^2 t}{\tau_t}$. 
        \\
        An associated Markov process is given by the solution of the SDE
            \begin{align*}
                \d Y_t = u_t \d t + \eta \d W_t;\quad Y_0 \sim \NN(x_0,\rho^2 I_d).
            \end{align*}
        \item $\mathcal{L}^{{\rm Jump}}_t f(x) = \int \big(f(y)  - f(x) \big) q_t (y,x) \d y$, where $q_t^{x_0,x_1}(y,x)
        = q_t(y,x) = \coloneqq \lambda_t(x) J_t(y)$,
        and
        \begin{align}
             J_t
             &\coloneqq \frac{\max(0,\xi_t) p_t}{\int \max(0,\xi_t) p_t \d y},
             \quad \lambda_t \coloneqq  \max\left(0,-\xi_t \right),
             \\
            \xi_t(y)
            &\coloneqq
            \frac{\eta^2(1-2t)}{2\tau_t}
            \Big(
            \frac{\|y-m_t\|^2}{\tau_t}  + 2 a_t^\tT \frac{y-m_t}{\sqrt{\tau_t}} - d
            \Big), \quad
            a_t \coloneqq \frac{\sqrt{\tau_t}}{\eta^2(1-2t)}(x_1-x_0).
        \end{align}
        Note that $\xi_t$ is well-defined for all $t \in [0,1]$ (the apparent singularity in $a_t$ is canceled by the prefactor) and that $\int \max(0,\xi_t) p_t \d y>0$ unless $t=1/2$ and $x_0=x_1$. In that case, $\lambda_t=0$ everywhere, and we set $J_{1/2}\coloneqq \mathcal{N}(0,I_d)$ for completeness.
        \item  If $\rho = 0$, then $\mathcal{L}^{{\rm Diff}}_t$ in i) simplifies to $u_t(x) =\frac{x_1-x}{1-t}$ and the SDE becomes
        \begin{align*}
            \d Y_t = \frac{x_1-Y_t}{1-t} \d t + \eta \d W_t, \quad Y_0=x_0.
        \end{align*}
    \end{romanenumerate}
\end{theorem}

\begin{figure}
  \centering
  \includegraphics[width=.21\textwidth]{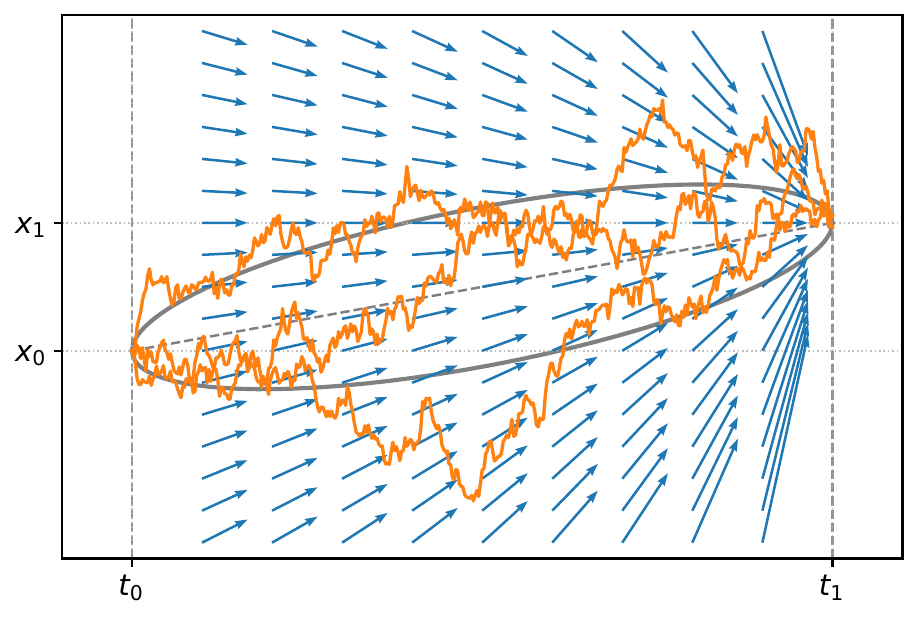}
  \hfill
  \includegraphics[width=.21\textwidth]{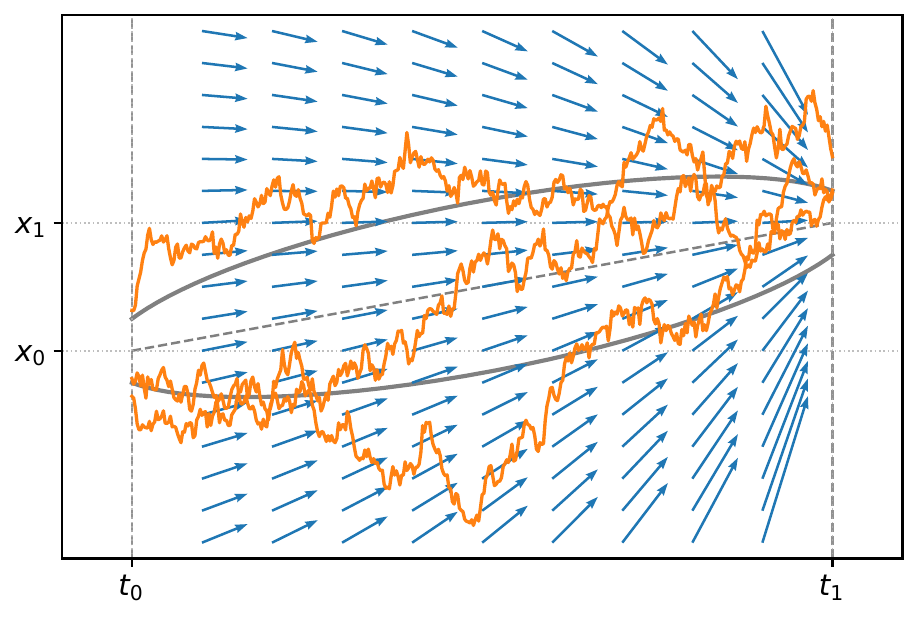}
  \hfill
    \includegraphics[width=.24\textwidth]{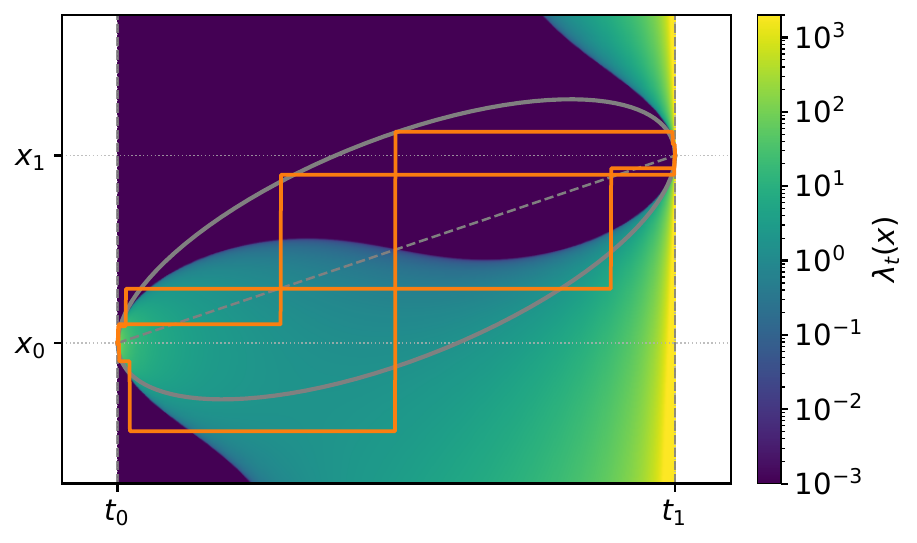}
  \hfill
  \includegraphics[width=.24\textwidth]{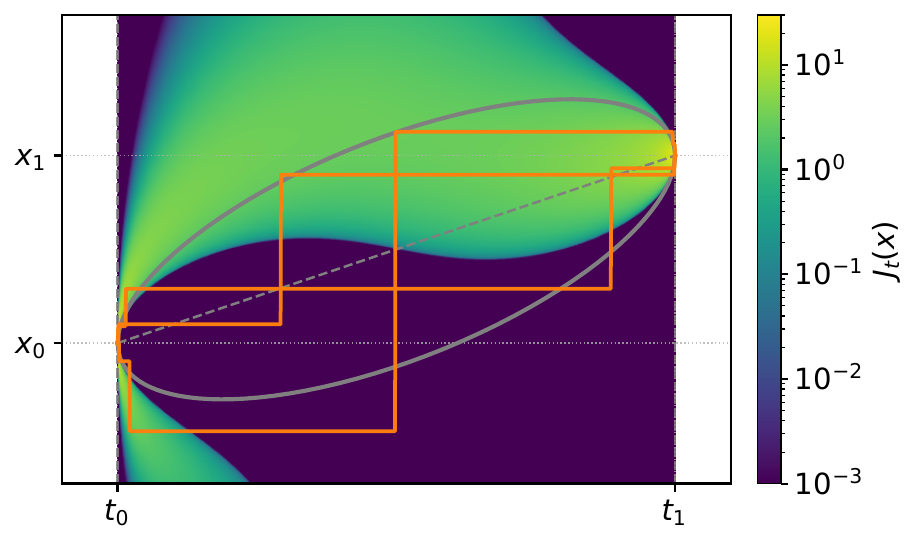}
  \caption{On the left side we see three paths of the interpolating diffusion bridge between $\mathcal{N}(x_0,\rho^2)$ and $\mathcal{N}(x_1,\rho^2)$, for small (leftmost) and large $\rho$ (left) with the underlying conditional drift (blue). It is evident that the latter drift field is less singular at the end point. On the right side we see three paths of the interpolating jump bridge together with heat maps of the jump intensity/jump probability $\lambda_t$ (right) and target/jump distribution $J_t$ (rightmost). In all cases we depict sample paths (orange) and marginal distribution $p_t$ (mean and standard deviation indicated by gray lines).}
  \label{fig:prop2}
\end{figure}

For $\rho=0$, we recover the drift of \cite{zhang2024trajectory}, which is singular at $t=1$. Moreover,  $P_0(\cdot,x_0,x_1) = \delta_{x_0}, P_1(\cdot,x_0,x_1) = \delta_{x_1},$
and we see that the Markov process generated by $\LL_t^{\rm Diff}$ is a Brownian bridge connecting $x_0$ and $x_1$. 
For $\rho > 0$, we don't encounter singularities and have
$P_0(\cdot,x_0,x_1) = \mathcal N(x_0,\rho^2 I_d)$, $P_1(\cdot,x_0,x_1) = \mathcal N(x_1,\rho^2 I_d).$ The diffusion and jump bridges are illustrated in Figure \ref{fig:prop2}.

\subsection{Loss Functions}
Let $\mathcal L^{{\rm Diff},x_0}_t$ and $\mathcal L^{{\rm Jump},x_0}_t$ denote the aggregated generators defined via \eqref{int} with respect to the conditional generators $\mathcal L_t^{{\rm Diff},x_0,x_1}$ and $\mathcal L_t^{{\rm Jump},x_0,x_1}$ from Theorem \ref{th1}. In this subsection, we derive a tractable loss function that allows to train neural parametrizations of these generators. 

As a first step, based on \cite[Proposition 2 and 3]{holderrieth2025generator}, we give a characterization of the aggregated generators as minimizers.
\begin{proposition} \label{prop:loss} 
    For $x_0, x_1 \in \R^d$ let $u_t^{x_0,x_1}$ and $q_t^{x_0,x_1}$ be the conditional drift and jump kernels from Theorem \ref{th1}. Define
    \begin{align}\label{loss_u}
    	 u_t^{x_0} \vcentcolon = \argmin_{v_t^{x_0} \in L_2(P_t^{x_0})} 
        \E_{t\in[0,1],x_1\sim P_{X_1|X_0=x_0}, x \sim P_{t}(\cdot,x_0,x_1)}
        \left[\| u_t^{x_0,x_1} (x) - v_t^{x_0}(x)\|^2
        \right]
    \end{align}
    and
    \begin{align}\label{loss_q}
    	 q_t^{x_0} \vcentcolon= \argmin_{r_t^{x_0}} 
        \E_{t\in[0,1],x_1\sim P_{X_1|X_0=x_0}, x \sim P_{t}(\cdot,x_0,x_1)}
        \left[\KL \big( q_t^{x_0,x_1}(\cdot,x), r_t^{x_0}(\cdot,x)\big)
        \right]
    \end{align}
    with the Kullback-Leibler divergence $\KL$, and the Markov kernel $P_t(\cdot,\cdot,\cdot)$ defined in Remark \ref{prop:post_x_0}. 
    Then, $\mathcal L^{{\rm Diff},x_0}_t$ and $\mathcal L^{{\rm Jump},x_0}_t$ are given by $(u_t^{x_0},\eta,0)$ and $(0,0,q_t^{x_0})$, respectively. Moreover, for each $\alpha\in[0,1]$ the process defined by the convex combination of the generators  
    \begin{align*}
        \LL_{t,\alpha}^{x_0}f\coloneqq\alpha \mathcal L^{{\rm Diff},x_0}_t f+ (1-\alpha)\mathcal L^{{\rm Jump},x_0}_t f
    \end{align*}
    has the marginal distribution $P_t^{x_0}$.
\end{proposition}

Note that we have samples from $P_{X_0,X_1}$ and not from $P_{X_1|X_0=x_0}$. So we integrate \eqref{loss_u} resp. \eqref{loss_q} with respect to $P_{X_0}$ and obtain the tractable loss functions
\begin{align}\label{eq:loss_general_generator}
    \mathcal E^{{\rm Diff}}(\theta) 
    &\coloneqq \E_{t\sim [0,1],(x_0,x_1)\sim P_{X_0, X_1}, x \sim P_{t}(\cdot,x_0,x_1)}
        \left[\| u_t^{x_0,x_1} (x)- v_t^{x_0,\theta}(x)\|^2 \right]\\\notag
    &=\E_{t\sim [0,1],(x,x_0,x_1)\sim P_t(\cdot,x_0,x_1)\times_{x_0,x_1}P_{X_0,X_1}}
        \left[ \| u_t^{x_0,x_1} (x)- v_t^{x_0,\theta}(x)\|^2 \right]\\
    \mathcal E^{{\rm Jump}}(\theta)
    &=\E_{t\sim [0,1],(x,x_0,x_1)\sim P_t(\cdot,x_0,x_1)\times_{x_0,x_1}P_{X_0,X_1}} \left[ \KL\big( q_t^{x_0,x_1}(\cdot,x), r_t^{x_0,\theta}(\cdot,x)\big) \right], \label{eq:loss_jump}
\end{align}
for neural parameterizations $v_t^{x_0,\theta}$ and $r_t^{x_0,\theta}$ of the velocity field and the jump kernel. Note that the condition $x_0$ is an additional input of the neural networks.

While learning the vector field $v_t^{\cdot, \theta}$ is straightforward, it is non trivial to choose a suitable parameterization of the rate kernel $r_t^{\cdot,\theta}$ and to compute the KL-divergence efficiently. The authors in \cite[Appendix F]{holderrieth2025generator} suggest to approximate the function  $r_t^{x_0,\theta}$ by dividing $\R$ into $n_{b}$ bins which results in an output dimension of $n_{b}\cdot d$ on $\R^d$ and is computationally demanding.
In this paper, we circumvent this difficulty by restricting our attention to functions $r_t^{x_0,\theta}$ of the form $r_t^{x_0,\theta}(\cdot, x)= \lambda_t^{x_0,\theta}(x) \mathcal N(\mu_t^{x_0,\theta}(x), \Sigma_t^{x_0,\theta}(x))$.
By the following theorem, the KL-divergence between $q_t^{x_0,x_1}(\cdot,x)$ and such functions can be computed explicitly.

\begin{theorem}\label{prop2}
    Let $q_t^{x_0,x_1}(\cdot,x)=\lambda_t(x)J_t(\cdot)$ be defined by Theorem \ref{th1}ii) with $\rho >0$, and denote by $\mu^{J}$ and $\Sigma^{J}$ the mean and covariance matrix of $J_t$.
    Then it holds for fixed $t,x \in [0,1] \times \R^d$ that
    \begin{align*} \label{kl_simple}
        &\KL(q_t^{x_0,x_1}(\cdot,x),\lambda \mathcal N(\mu, \Sigma)) 
        = 
        F_{t,x} (\lambda,\mu,\Sigma) + {\rm const},
        \\\notag
         &F_{t,x} (\lambda,\mu,\Sigma) 
         \coloneqq 
         \lambda + \frac{\lambda_t(x)}{2} \left( \log(\det\Sigma) - 2\log(\lambda) 
        + \left\langle\Sigma^{-1},\Sigma^J\right\rangle + \left\langle \mu^J-\mu,\Sigma^{-1}(\mu^J-\mu)\right\rangle\right),
    \end{align*}
    where ${\rm const}$ is a constant independent of 
    $(\lambda,\mu,\Sigma)$. We note that, for $\lambda_t(x)\neq0$ and positive definite $\Sigma^J$, the problem
    \[
    \min_{(\lambda,\mu,\Sigma)\in\R_{>0} \times \R^d \times \textup{SPD}(d)}
    F_{t,x} (\lambda,\mu,\Sigma) 
    \]
    has the unique minimizer $(\lambda_t(x), \mu^{J}, \Sigma^{J})$.
\end{theorem}

Combining \eqref{eq:loss_jump} and Theorem \ref{prop2}, we learn the jump process, which allows us to approximately sample from $P_{X_1|X_0=x_0}$, by minimizing the loss 

\begin{equation} \label{eq:loss_J}
    \mathcal E^{{\rm Jump}} (\theta) 
    \coloneqq \E_{t\sim [0,1],(x,x_0,x_1)\sim P_t(\cdot,x_0,x_1)\times_{x_0,x_1}P_{X_0,X_1}}\Big[ F_{t,x}\big(\lambda_t^{x_0,\theta}(x),\mu_t^{x_0,\theta}(x),\Sigma_t^{x_0,\theta}(x)\big) \Big].
\end{equation}

In order to evaluate $F_{t,x}$,
we need to compute the mean and covariance matrix $(\mu^{J},\Sigma^{J})$ of $J_t$. Surprisingly, we can derive closed-form formulas for these in dimension $d\in\{1,2\}$.

\begin{theorem}\label{prop1}
    For dimensions $d\in\{1,2\}$ the mean $\mu^{ J}=\mu_{x_0,x_1}^J$ and variance $\Sigma^J=\Sigma^J_{x_0,x_1}$ of $J_t$ can be computed in closed-form, see Theorem \ref{prop:formula1d} and Theorem \ref{prop:2d_formulas}.
\end{theorem}

In dimension $d\geq 3$, formulas for $\mu^J$ and $\Sigma^J$ can still be obtained, so that training according to \eqref{eq:loss_J} would in principle be possible. But there is an alternative approach, which was numerically superior in our experiments in Section \ref{sec:numerics}.

\begin{remark}[Factorized jump paths.]\label{rem:factorized}
    Instead of computing $\lambda_t,\mu^J,\Sigma^J$ in higher dimensions, we can also use the one dimensional objects and stack them, see \cite[Proposition 5]{holderrieth2025generator}. More precisely, write $x_0=(x_0^1,\ldots,x_0^d),\, x_1=(x_1^1,\ldots,x_1^d)$ and let $P_t^i=P_t(\cdot,x_0^i,x_1^i)$ be the one dimensional paths with corresponding jump kernel $q_t^{x_0^i,x_1^i}$ from Theorem \ref{th1}. Let $P_t=\prod_{i=1}^d P_t^i$, then $P_0=\NN(x_0,\rho^2 I_d)$ and $P_1=\NN(x_1,\rho^2 I_d)$. Now take  
    \begin{align}
    	q_t^{i,x_0} = \argmin_{r_t^{x_0}} 
        \E_{t\in[0,1],x_1\sim P_{X_1|X_0=x_0}, x \sim P_{t}(\cdot,x_0,x_1)}
        \left[\KL \big( q_t^{x_0^i,x_1^i}(\cdot,x^i), r_t^{x_0}(\cdot,x)\big) \right],
    \end{align}
    where minimization is over all jump kernels $r_t^{x_0}: \mathcal{B}(\R) \times \R^d \to [0,\infty)$. Then $q^{x_0}_t\coloneqq (q^{1,x_0}_t, \ldots,q^{d,x_0}_t)$ also determines a process $Y_t^{x_0}$ with $Y_0^{x_0}\sim\NN(x_0,\rho^2 I_d)$ and $Y_1^{x_0}\sim P_{X_1|X_0=x_0}*\NN(0,\rho^2 I_d)$. Sampling from this process can be done by sampling each coordinate for each time step independently as in Algorithm \ref{alg:markov_from_generator}. Note that  $q^{i,x_0}_t$ depends on the whole state $x\in\R^d$. Therefore, after each update the input of $q_t^{i,x_0}$ changes in every coordinate, and the induced process is not just the product of one dimensional processes.
\end{remark}

\section{Learning Irregular Time Series}\label{section:whole_times_series}
Now we formalize the procedure of learning entire time series based on the interpolants from Section \ref{sec:interpolating_mps} and the scheme outline in Section \ref{sec:problem_setting}, see also Figure \ref{fig:conditional_generator}. 

We train two networks $v_{\cdot}^{\cdot,\cdot,\theta}$ and $r_{\cdot}^{\cdot,\cdot,\theta}$ that approximate the drift and jump kernel of interpolating processes on arbitrary time intervals. Inputs of these networks are the current position and time, the end time of the current interval, and the memory. 

For a time series $x^{i}\in \mathcal{X}$ let $x_j^i \coloneqq x^i_{t^i_j}$ and define the memory
$\xi^i_j \coloneqq \big( (x^i_{0},t^i_0), \ldots, (x^i_j,t^i_j)\big)$.
For $t\in[0,T]$, there exists $j \in\{0,\ldots,n_i-1\}$ such that $t \in [ t^i_j, t^i_{j+1})$. We denote this dependence by $j=j(t,i)$.  Then we learn the network for the drift
by minimizing the loss function
\begin{align*}
    \mathcal{E}^{\rm Diff}(\theta):= \E_{t \sim \mathcal{U}\left([0,T
    ]\right), i \sim \mathcal{U}\left([N]\right), x\sim P_t\left(\cdot,x^i_{j}, x^{i}_{j+1}\right)}
    \Big[\big\| u_t^{x^{i}_{j},x^{i}_{j+1}}(x)- v_t^{\xi_j^i,t^i_{j+1},\theta}(x)\big\|^2\Big], 
\end{align*}
where $u_t^{x_j^i,x_{j+1}^i}$ is the generator of the Markov process connecting $x_{j(t,i)}^i$ and $x_{j(t,i)+1}^i$ from Theorem \ref{th1}i). Similarly, we set up a network to learn the jump rate kernel with loss
\begin{align}
    \mathcal{E}^{\rm Jump}(\theta):= \E_{t \sim \mathcal{U}\left([0,T
    ]\right), i \sim \mathcal{U}\left([N]\right), x\sim P_t\left(\cdot,x^i_{j}, x^{i}_{j+1}\right)}\Big[{\rm KL}\big( q_t^{x^{i}_{j},x^{i}_{j+1}}(x),r_t^{\xi_j^i,t^i_{j+1},\theta}(x)\big) \Big], 
\end{align}
where $q_t^{x_j^i,x_{j+1}^i}$ is the generator of the Markov process connecting $x_j^i$ and $x_{j+1}^i$ from Theorem \ref{th1}ii). (The formulation of Theorem \ref{th1} for general time intervals $[t_j,t_{j+1}]$ is given in Appendix \ref{paragraph:th1}.)
As explained above, we parametrize Gaussian-like jump measures and compute the loss in closed form using Theorems \ref{prop2} and \ref{prop1}.

After learning the velocity fields and jump kernels, we sample our time series according to Algorithm \ref{alg:glueing_markov_proccesses}.
This is based on Algorithm \ref{alg:markov_from_generator} with time
intervals $[t_i,t_{i+1}]$, $i=0, \ldots, n-1$ instead of $[0,1]$. 
Here, the generator $\LL_t^{\cdot,\cdot}$ is given by the triplet $\bigl(\alpha \,v_t^{\cdot,\cdot,\theta},\sqrt{\alpha}\,\eta, (1-\alpha)\, r_t^{\cdot,\cdot,\theta}\bigr)$ for $\alpha\in[0,1]$, where $\alpha=0$ corresponds to pure jump interpolants, $\alpha=1$ corresponds to drift-diffusion interpolants, and $\alpha\in (0,1)$ to a mixture (Markov superposition).

\begin{algorithm}[H]
\caption{Sampling from glued Markov processes at $\tau=(t_0, \ldots,t_n)$.}
\label{alg:glueing_markov_proccesses}
\begin{algorithmic}[1]
    \State \textbf{Given:} $\LL_t^{\cdot,\cdot}$ and $x_0 \sim X_0$, $\tau=(t_0, \ldots, t_n)$. $\text{ALG}_1$ from Algorithm~\ref{alg:markov_from_generator}.
    \State \textbf{Initialize:} $\xi_0=(x_0,t_0)$
    
    \For{$i = 0$ \textbf{to} $n-1$}
        \State $x_{i+1} = \text{ALG}_1(x_i, \LL_t^{\xi_i,t_{i+1}})$
        \State $\xi_{i+1} = (\xi_i, (x_{i+1}, t_{i+1}))$
    \EndFor
    
    \Return $(x_0,\ldots,x_n)$ approximately distributed as $(X_{t_0},\ldots,X_{t_n})$
\end{algorithmic}
\end{algorithm}

For this autoregressive sampling procedure, the approximation error for the finite-dimensional marginal distributions can be controlled according to the following stability result. A precise statement and proof are given in Appendix~\ref{app:joint_post}.

\begin{theorem}[Stability of finite-dimensional marginals - informal]
\label{rem:joint_post}
Let $\mu_i$ be a Markov kernel approximating $P_{X_{i+1}|X_0,\ldots ,X_i}$ for $i=0,\ldots,n-1$. Set $ \hat{P} \coloneqq \mu_{n-1}(\cdot,x_0,\ldots,x_{n-1})\times_{x_0,\ldots,x_{n-1}}\cdots\times_{x_0}P_{0}\in \PP_2(\R^{(n+1)d})$ for some $P_0\in \PP_2(\R^{d})$. If each $\mu_i$ is Lipschitz-continuous in $(x_0,\ldots, x_i)$ with respect to the Wasserstein-2 metric, then the approximation error $W_2(P_{X_0,\ldots,X_n},\hat{P})$ can be controlled in terms of the initial approximation error $W_2(P_{X_0},P_0)$ and the approximation error between the conditional laws and Markov kernels.
\end{theorem}

\section{Numerical Illustrations}\label{sec:numerics}
We evaluate the proposed approach on synthetic time series in low and high dimension. The low-dimensional experiments use one- and two-dimensional Black-Scholes stock data,  the high-dimensional experiment uses a moving-box video dataset, designed to exhibit both regular continuous motion and abrupt jumps.

In the experiments we do not use the entire process memory. Instead, we fix a memory length $m$ and feed only the truncated memory $\xi=\bigl((x_{j-m},t_{j-m}),\ldots,(x_j,t_j)\bigr)$ into the networks, where we set $(x_k,t_k)=(x_0,0)$ for $k<0$. For the low-dimensional experiments, the drift and jump networks are implemented as simple feedforward neural networks and trained with the Adam optimizer \cite{adam}. For the video experiment, we use U-Net based architectures.

\subsection{Synthetic Black-Scholes Data}
\label{sec:bs_data}

\begin{figure}
    \centering
    
    \includegraphics[width = 0.24\textwidth, trim=15mm 5mm 15mm 5mm,
    clip]{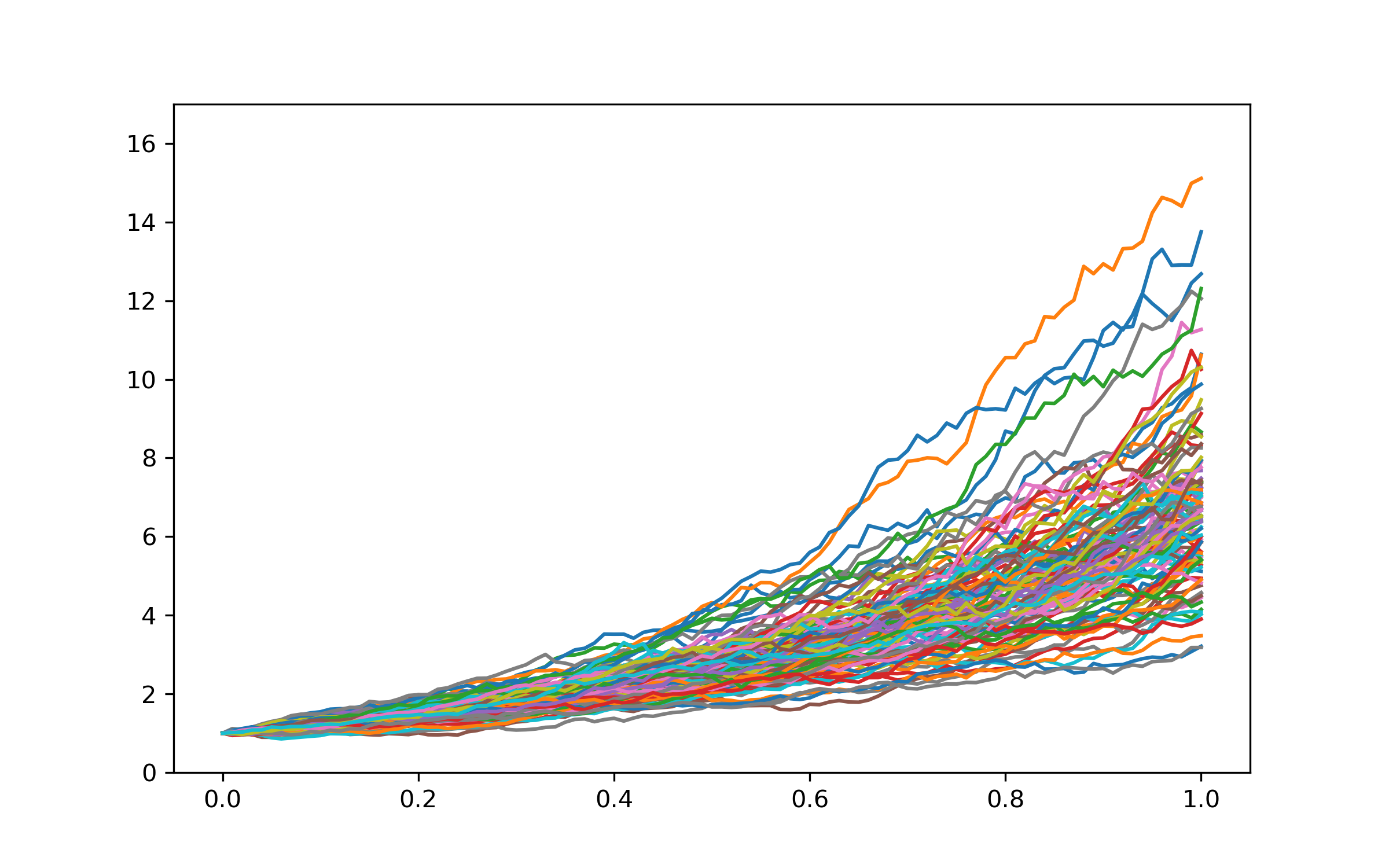}
    \includegraphics[width = 0.24\textwidth, trim=15mm 5mm 15mm 5mm,
    clip]{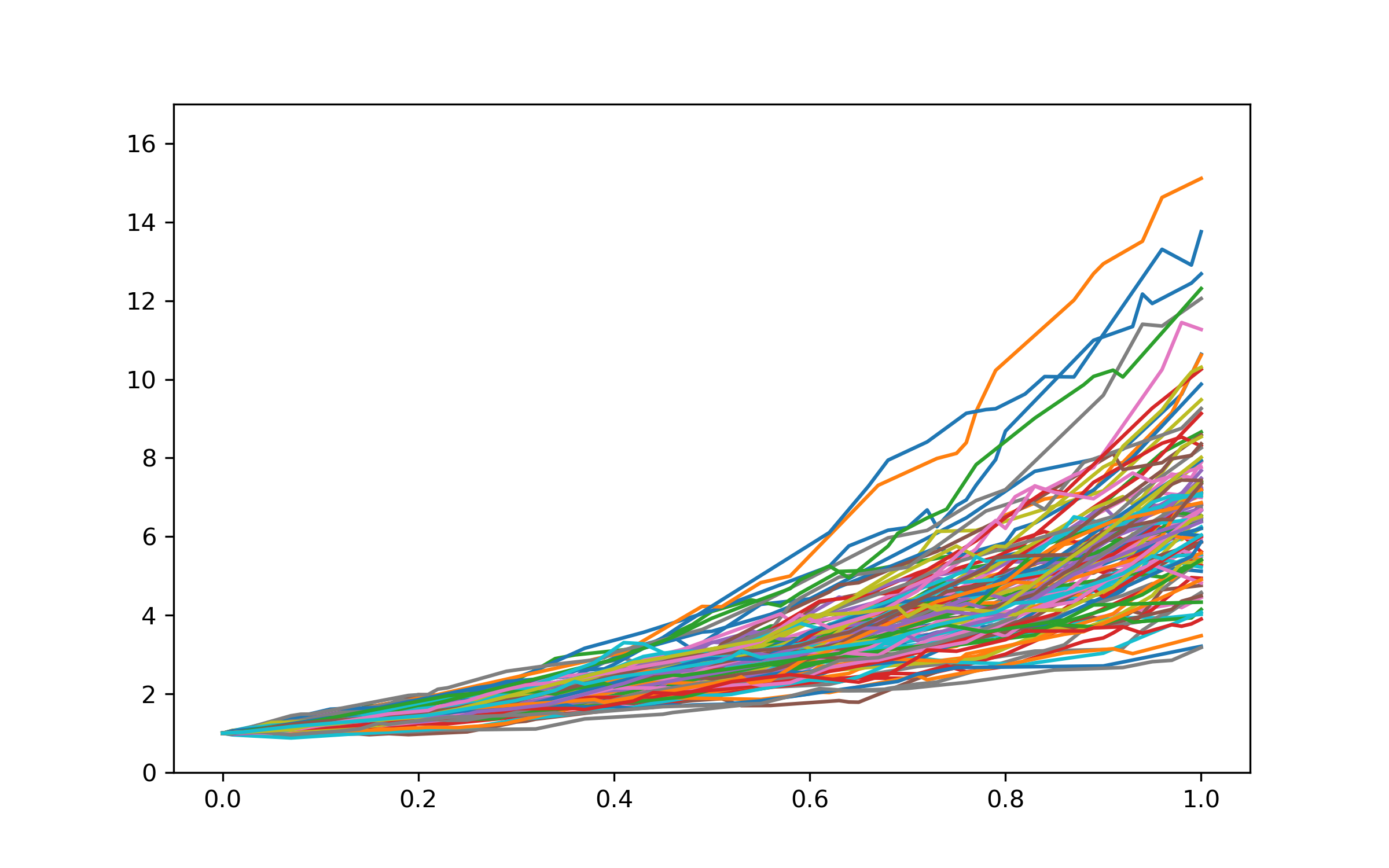}
        \includegraphics[width = 0.24\textwidth, trim=15mm 5mm 15mm 5mm,
    clip]{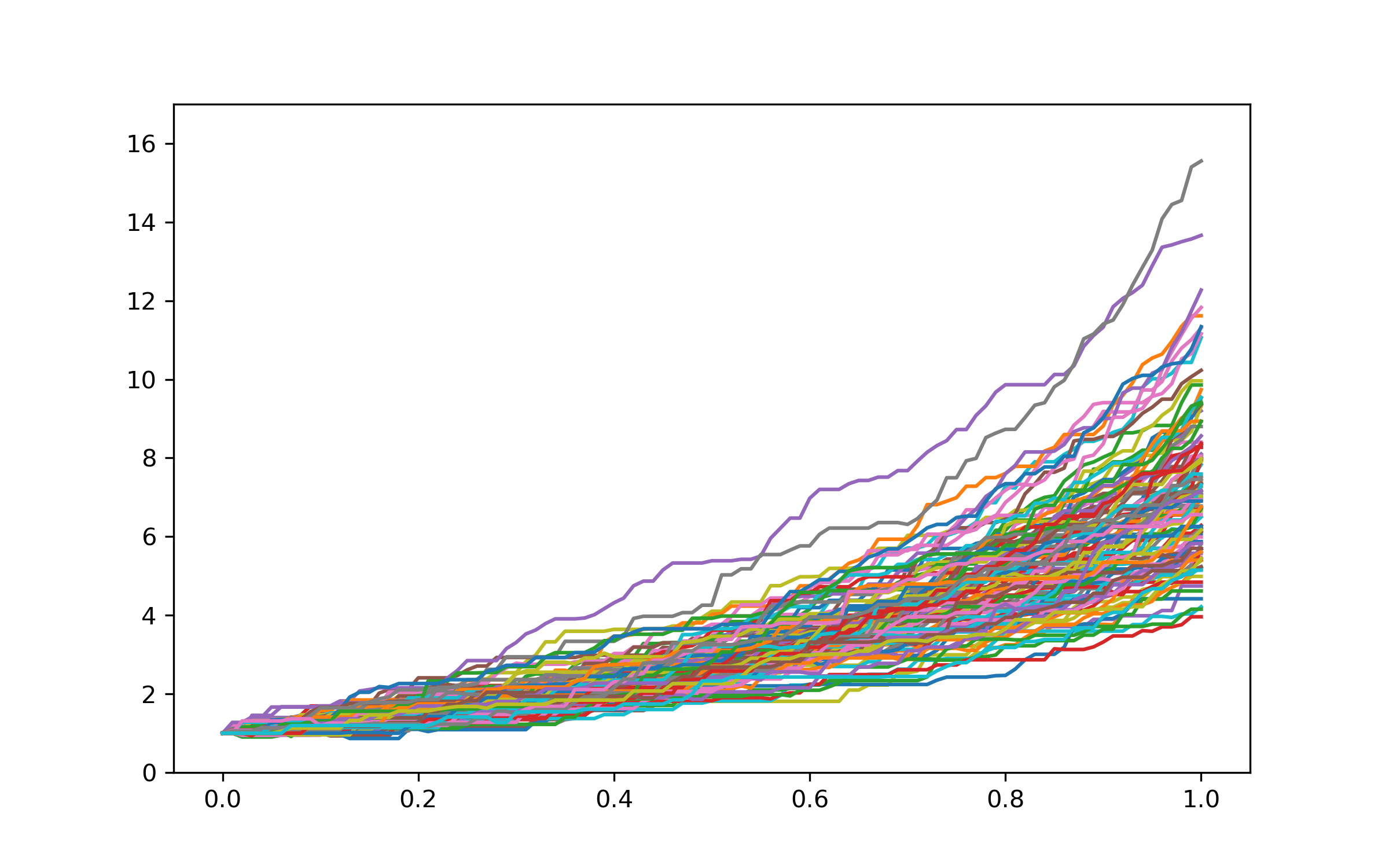}
    \par
    \includegraphics[width = 0.24\textwidth, trim=15mm 5mm 15mm 5mm,
    clip]{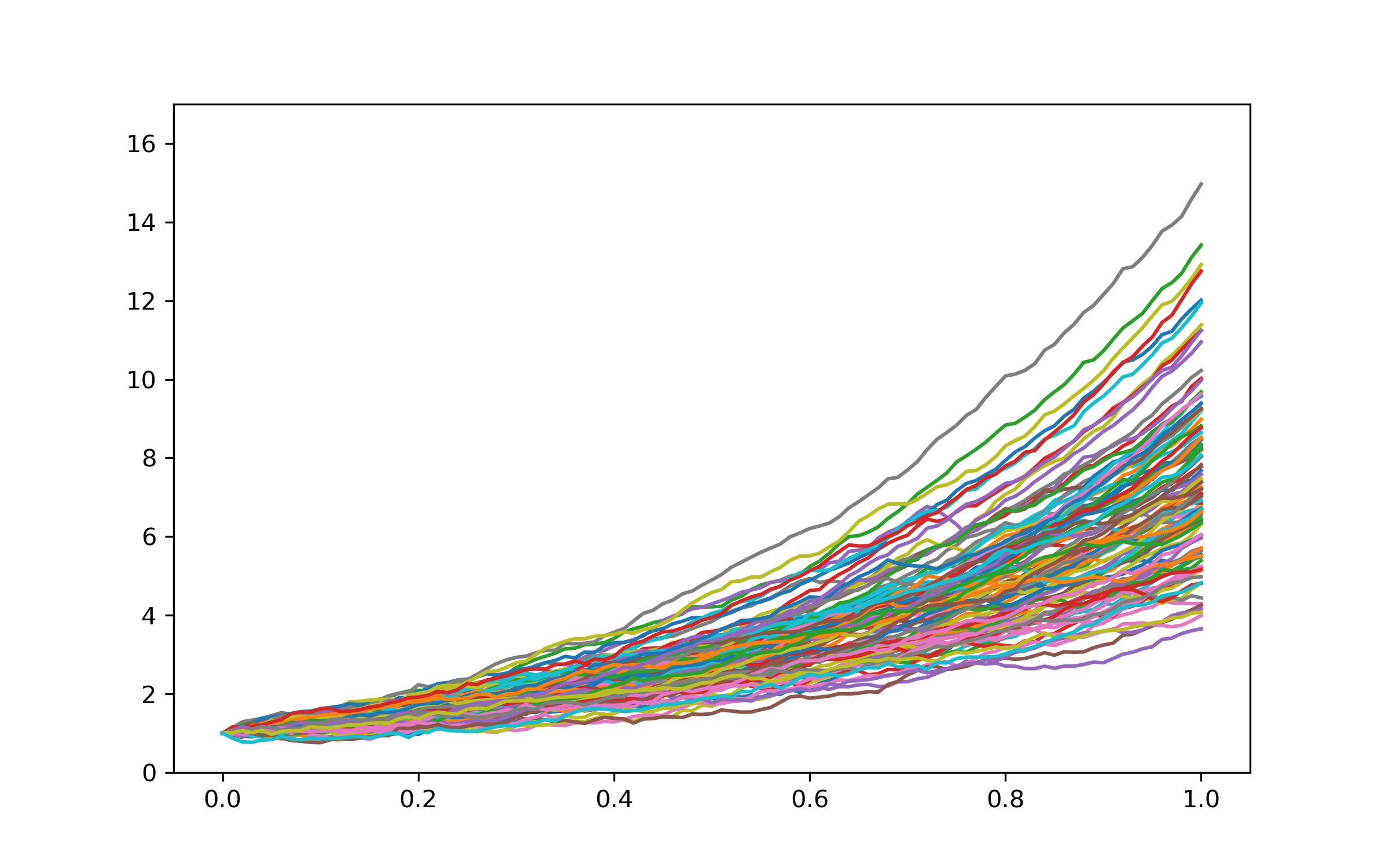}
    \includegraphics[width = 0.24\textwidth, trim=15mm 5mm 15mm 5mm,
    clip]{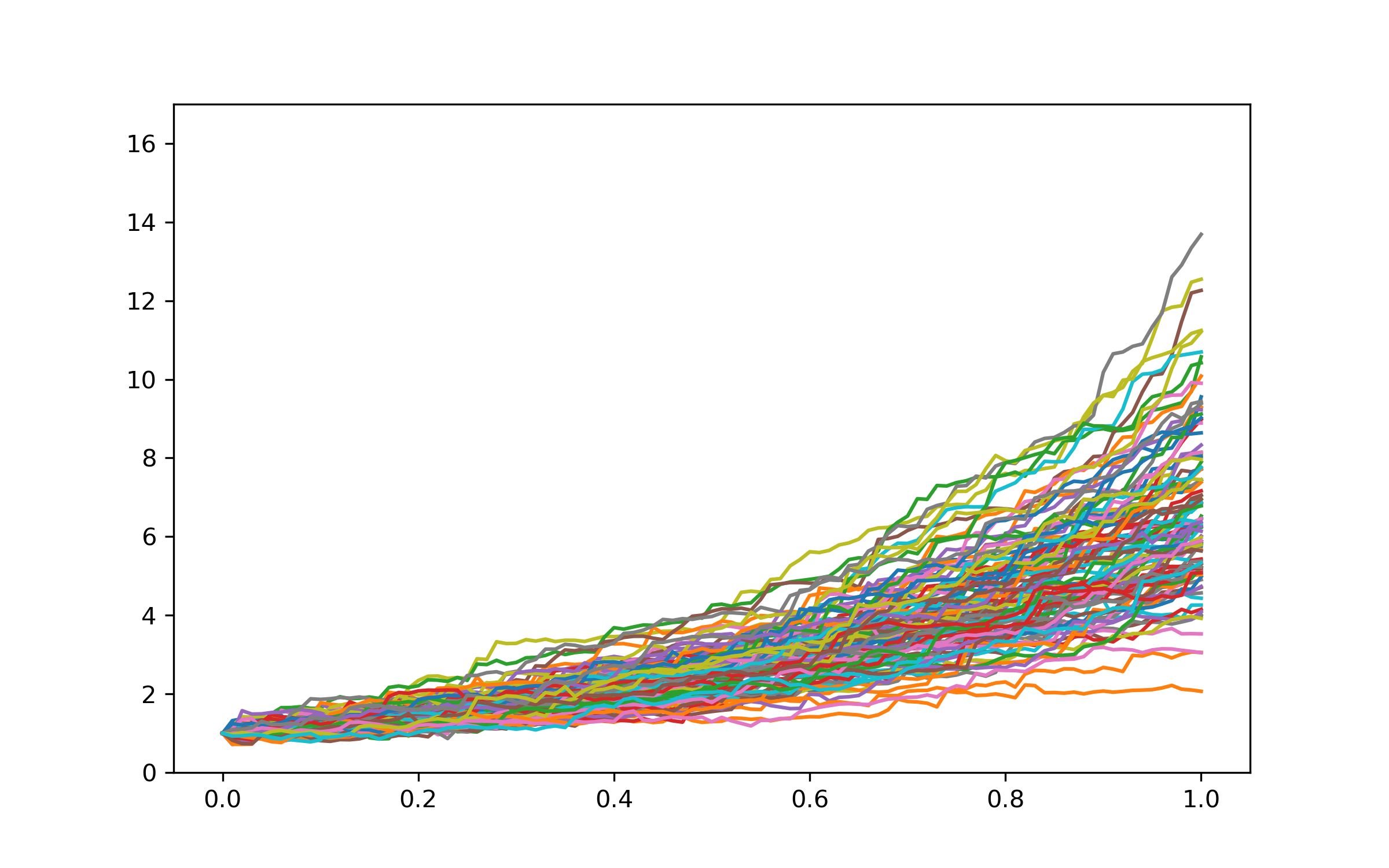}
    \includegraphics[width = 0.24\textwidth, trim=15mm 5mm 15mm 5mm,
    clip]{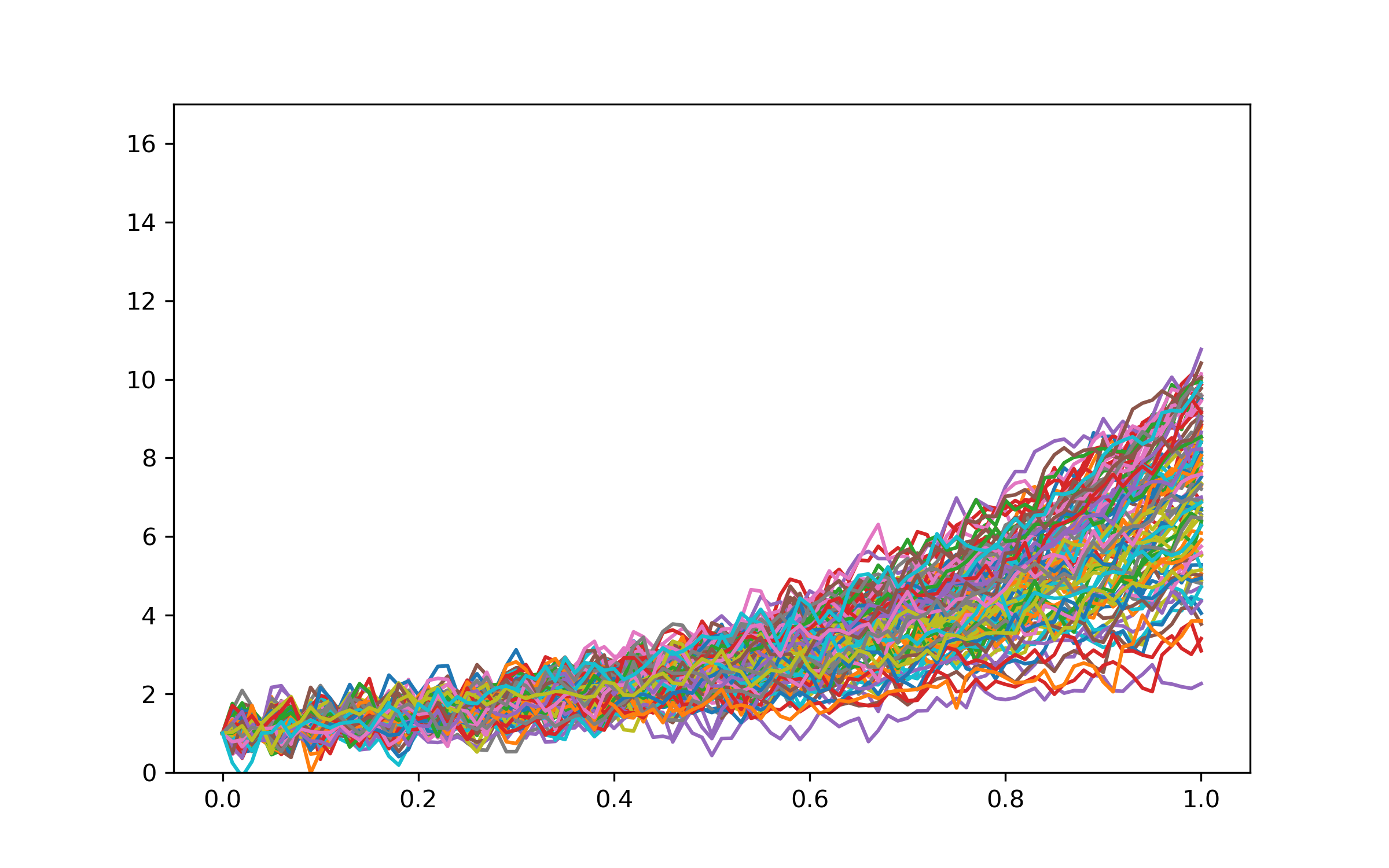}
    \includegraphics[width = 0.24\textwidth, trim=15mm 5mm 15mm 5mm,
    clip]{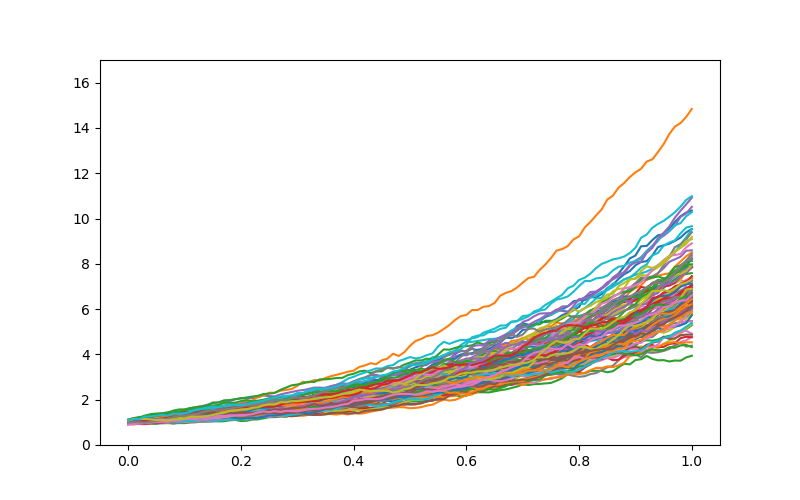}
    \caption{Illustration of '25\%'-column in Table \ref{tab:mmd_subsampling:stock}. From left to right, top to bottom: Ground truth trajectories of one-dimensional Black-Scholes dataset. Irregular, non-aligned trajectories after 25\% subsampling. Trajectories generated by the following models trained on data in the 25\% observation regime: JUMP, SDE,  MS ($\alpha=0.2$), TFM, SDE matching.}
    \label{fig:Black-Scholes-sub25}
\end{figure}

We use a Black-Scholes model with fixed parameters from \cite{herrera2021neural} to generate synthetic stock data on the uniform time grid $\{t=\tfrac{i}{100} : i=0,\ldots, 100\}$ in one and two dimensions. Visualizations are given in Figure \ref{fig:Black-Scholes-sub25} and the supplementary material.

To study irregular and non-aligned observation regimes, we randomly subsample these trajectories while always retaining the initial and terminal time points. We consider trajectories observed at 11, 26, and 101 time points, corresponding to 10, 25, and 100 time intervals of varying length. 

During generation according to Algorithm \ref{alg:glueing_markov_proccesses}, we choose 11, 26, or 101 equidistant time points, i.e., 10, 25 or 100 regular bridges, respectively, and in Algorithm \ref{alg:markov_from_generator} the step size $h=0.001$. The resulting grid aligns with the original 101-point grid, so that all generated trajectories can be compared to fully observed ground truth paths.

As a metric we report the Maximum Mean Discrepancy (MMD) with negative distance kernel \cite{szekely} between the generated trajectories and ground truth (test) trajectories using the geomloss package \cite{feydy2019interpolating}. The same full-trajectory MMD between generated and ground truth (validation) trajectories is used for checkpoint selection and hyperparameter tuning.

We compare our jump-based method, our SDE-based method and their Markov superposition with the TFM method from \cite{zhang2024trajectory} and SDE matching \cite{bartosh2025sde}. 
Further possible baselines would be (latent) NeuralODE/SDE \cite{chen2018,rubanova2019latent, li2020scalable, kidger2021efficient} or aligned Flow Matching \cite{liu2023i2sb, somnath2023aligned}. However we did not consider them explicitly here, as those methods were consistently outperformed by TFM in experiments in \cite{zhang2024trajectory}.

For SDE matching, we modified the publicly available code to make it applicable to irregular observation times. For the TFM baseline, we used the public repository of~\cite{zhang2024trajectory} as a reference. In preliminary experiments, we found that the released implementation combines several repository-specific modeling, training, and rollout choices that differ from the paper-level formulation. Direct use of this implementation led to qualitatively implausible trajectories on the Black--Scholes benchmark. We therefore implemented the TFM method directly in our framework, following the mathematical description in~\cite{zhang2024trajectory} as closely as possible. In particular, we used the denoising objective proposed in \cite{zhang2024trajectory}, whereas for our SDE method we predict the velocity.

\begin{table}
\centering
\begin{minipage}{\linewidth}
\centering
\scalebox{0.65}{
\begin{tabular}{llll}
\hline
\textbf{Method} & \textbf{10\% observed data} & \textbf{25\% observed data} & \textbf{100\% observed data} \\
\hline
      Jump-based method
      & [0.011 {\scriptsize$\pm$ 0.0035}] {\scriptsize$(\eta^2=0.3)$}
      & [\textbf{0.012} {\scriptsize$\pm$ 0.0059}] {\scriptsize$(\eta^2=0.3)$}
      & [\textbf{0.015} {\scriptsize$\pm$ 0.0057}] {\scriptsize$(\eta^2=0.3)$} \\

      SDE-based method
      & [\textbf{0.0076} {\scriptsize$\pm$ 0.0015}] {\scriptsize$(\eta^2=0.3)$}
      & [0.022 {\scriptsize$\pm$ 0.0016}] {\scriptsize$(\eta^2=0.3)$}
      & [0.027 {\scriptsize$\pm$ 0.0061}] {\scriptsize$(\eta^2=0.3)$} \\

      Jump + SDE (Markov superposition)
        & [\textbf{0.0065} {\scriptsize$\pm$ 0.0030}] {\scriptsize$(\eta^2=0.3,\alpha=0.9)$}
        & [\textbf{0.010} {\scriptsize$\pm$ 0.0051}] {\scriptsize$(\eta^2=1.0,\alpha=0.2)$}
        & [\textbf{0.012} {\scriptsize$\pm$ 0.0069}] {\scriptsize$(\eta^2=0.3,\alpha=0.3)$} \\

      TFM method \cite{zhang2024trajectory}
      & [0.064 {\scriptsize$\pm$ 0.012}] {\scriptsize$(\eta^2=0.3)$}
      & [0.080 {\scriptsize$\pm$ 0.050}] {\scriptsize$(\eta^2=10.0)$}
      & [0.026 {\scriptsize$\pm$ 0.0084}] {\scriptsize$(\eta^2=1.0)$} \\

      SDE matching \cite{bartosh2025sde}
      & [0.084 {\scriptsize$\pm$ 0.016}]
      & [0.091 {\scriptsize$\pm$ 0.009}]
      & [0.029 {\scriptsize$\pm$ 0.008}] \\
      \hline
\end{tabular}
}
\vspace{1ex}
\caption{Test MMD over 5 runs on the one-dimensional Black-Scholes dataset. The selected model noise parameter is shown in parentheses; for Markov superpositions, $\alpha$ denotes the
SDE weight.}
\label{tab:mmd_subsampling:stock}
\end{minipage}
\end{table}

\begin{table}
\centering
\begin{minipage}{\linewidth}
\centering
\scalebox{0.59}{ 
\begin{tabular}{llll}
\hline
\textbf{Method} & \textbf{10\% observed data} & \textbf{25\% observed data} & \textbf{100\% observed data} \\
\hline

      Jump-based with general $\Sigma$
      & [0.198 {\scriptsize$\pm$ 0.072}] {\scriptsize$(\eta^2=3.0)$}
      & [0.101 {\scriptsize$\pm$ 0.017}] {\scriptsize$(\eta^2=3.0)$}
      & [0.282 {\scriptsize$\pm$ 0.076}] {\scriptsize$(\eta^2=0.3)$} \\

      Jump-based with factorized components
      & [0.032 {\scriptsize$\pm$ 0.0086}] {\scriptsize$(\eta^2=1.0)$}
      & [\textbf{0.026} {\scriptsize$\pm$ 0.0084}] {\scriptsize$(\eta^2=1.0)$}
      & [0.094 {\scriptsize$\pm$ 0.0036}] {\scriptsize$(\eta^2=0.3)$} \\

      SDE-based method
      & [0.023 {\scriptsize$\pm$ 0.0035}] {\scriptsize$(\eta^2=1.0)$}
      & [0.029 {\scriptsize$\pm$ 0.0087}] {\scriptsize$(\eta^2=1.0)$}
      & [\textbf{0.051} {\scriptsize$\pm$ 0.034}] {\scriptsize$(\eta^2=1.0)$} \\

      Jump (factorized) + SDE (Markov superposition)
        & [\textbf{0.021} {\scriptsize$\pm$ 0.0011}] {\scriptsize$(\eta^2=1.0,\alpha=0.8)$}
        & [\textbf{0.022} {\scriptsize$\pm$ 0.0064}] {\scriptsize$(\eta^2=1.0,\alpha=0.4)$}
        & [\textbf{0.046} {\scriptsize$\pm$ 0.035}] {\scriptsize$(\eta^2=1.0,\alpha=0.8)$} \\

      TFM method \cite{zhang2024trajectory}
      & [0.143 {\scriptsize$\pm$ 0.033}] {\scriptsize$(\eta^2=1.0)$}
      & [0.202 {\scriptsize$\pm$ 0.049}] {\scriptsize$(\eta^2=10.0)$}
      & [0.093 {\scriptsize$\pm$ 0.023}] {\scriptsize$(\eta^2=10.0)$} \\
      \hline
      
\end{tabular}
}
\vspace{1ex}
\caption{Test MMD over 5 runs on the two-dimensional Black-Scholes dataset. The selected model noise parameter is shown in parentheses; for Markov superpositions, $\alpha$ denotes the
SDE weight. The first method was trained according to \eqref{eq:loss_J} using Theorem \ref{prop:2d_formulas}, the second according to Remark \ref{rem:factorized}.}
\label{tab:mmd_subsampling:stock2d}
\end{minipage}
\end{table}

Tables \ref{tab:mmd_subsampling:stock} and \ref{tab:mmd_subsampling:stock2d} report the test MMD for different subsampling/observation regimes. In the subsampled regimes, our proposed jump and SDE methods substantially outperform the baselines. The Markov superposition often improves over both individual components, indicating the usefulness to combine complementary local mechanisms. In the fully observed regime, the gap to TFM and SDE matching narrows, but remains significant. 

Table \ref{tab:mmd_subsampling:stock2d} also shows that the jump-based method with factorized components, discussed in Remark~\ref{rem:factorized}, strongly outperforms the jump-based method using \eqref{eq:loss_J} and the explicit formulas for $\mu^J$ and $\Sigma^J$ in two dimensions.

It is evident from both tables, that performance of our proposed methods can decrease as the number of observation points increases. Two ablations are reported in the supplementary material investigating this issue. We attribute it to a large part to the chosen memory representation: Since Black-Scholes dynamics are Markovian, memory is not required in principle. A no-memory ablation for the one-dimensional case improves performance in the 25\%- and 100\% observation regimes, removing the observable trend almost entirely. This suggests that the flat memory representation used here may not be optimal when many observations are available.

Further details and plots are presented in Section \ref{sec:supp:stock_details} of the Supplementary Material.

\subsection{Moving-Box Video Data}
\label{sec:video_data}

To test our method on high-dimensional, non-Markovian data, we introduce a synthetic video dataset. Each state is a $16\times16$ binary image $x_t\in\{0,1\}^{256}$, showing a white $3x3$ box on black background. The box has a horizontal and a vertical direction of movement, which are randomly initialized along with its position. With each time step, the box moves horizontally by one pixel, and with probability $30\%$ it performs an additional vertical jump of three pixels. The horizontal movement is subject to reflections at the image boundary. The jump direction is reversed as soon as no complete jump of three pixels can be performed in the given direction.

We apply the jump-based method componentwise according to Remark \ref{rem:factorized}. During autoregressive generation, the memory is updated using a binarized version of the network output: the nine largest pixel values are set to one, the remaining ones are set to zero. This ensures that the model is evaluated in the binary image regime in which it was trained.

We compare our jump-based method, our SDE-based method, and superpositions thereof with the TFM method. We report the MMD with negative distance kernel between generated and ground truth (test) trajectories of the box's position. 

For suitable values of $\eta$, all models learn to respect the shape of the box and model its movement well. Lowest MMD is then achieved by the SDE-based model, followed by the jump-based model. TFM yields the highest MMD. The optimal Markov superposition has a high SDE-weight and performs similarly to the SDE-based model. Quantitative results and further details are given in Section \ref{sec:supp:details_video} of the Supplementary Material.

\section{Conclusions}
We introduced a new method to learn irregularly sampled time series. To this end, we constructed  stabilized bridges interpolating locally between subsequent observations. Those are based on specific diffusion and jump  processes, admitting explicit conditional generators suitable for simulation-free training. We proposed a Gaussian approximation for the jump kernel and derived closed-form formulas for the corresponding KL divergence. Our method uses a memory to recover the correct global distributions for non-Markovian dynamics. We tested its performance on datasets in various dimensions. 

One limitation of the current formulation is the restriction of the rate kernel to Gaussians. A next step would be to employ more complex models, such as k-parametric exponential families or Gaussian mixtures. 
Another promising extension would be to replace the empirical tuning of auxiliary parameters, such as the Markov superposition weight $\alpha$ or the diffusion $\eta$, by a principled approach. 
Moreover, our experiments indicate that the method could benefit from improved memory encodings, for example recurrent or attention-based representations. 

\appendix
\renewcommand*{\theHsection}{appendix.\Alph{section}}
\renewcommand*{\theHsubsection}{\theHsection.\arabic{subsection}}

\section{Proofs from Section \ref{sec:interpolating_mps} and General Interpolation Times} \label{app:proofs}

The results from Section \ref{sec:interpolating_mps} extend to general interpolation times $t_0$ and $t_1$ (not necessarily equal to 0 and 1) using the following definitions for $m_t$, $\tau_t$, $u_t$ and $\xi_t$:
\begin{itemize}
\item $m_t \coloneqq \frac{t_1-t}{t_1-t_0}x_0 + \frac{t-t_0}{t_1-t_0} x_1$ and
$\tau_t\coloneqq \eta^2 \frac{(t-t_0)(t_1-t)}{(t_1-t_0)}+\rho^2$,
\item $u_t \coloneqq \frac{x_1-x_0}{t_1-t_0} +  \left(x-m_t\right) \frac{\eta^2}{2\tau_t}\left(\frac{t_0+t_1-2t}{(t_1-t_0)}
-1\right)$,
\item $\xi_t(y) \coloneqq \frac{\eta^2(t_0+t_1-2t)}{2(t_1-t_0)\tau_t}\left(\frac{\|y-m_t\|^2}{\tau_t}+2 a_t^\tT  \frac{y-m_t}{\sqrt{\tau_t}} - d\right)$ with $a_t \coloneqq \frac{\sqrt{\tau_t}}{\eta^2(t_0+t_1-2t)}(x_1-x_0)$.
\end{itemize}
The definitions of $P_t$, $p_t$, $J_t$ and $\lambda_t$ in terms of the quantities above remain unchanged, except that $J_t \coloneqq \mathcal{N}(0,I_d)$ for $t=(t_0+t_1)/2$ instead of $t=1/2$.
Below we prove Theorem \ref{th1}, Theorem \ref{prop2} and Theorem \ref{prop1} in this general setting.

\begin{proof}[Proof of Theorem \ref{th1}.]\label{paragraph:th1}
By the definitions above, we have that $p_{t_0} \sim \mathcal N(x_0,\rho^2 I_d)$ and $p_{t_1} \sim \mathcal N(x_1, \rho^2 I_d)$. Moreover, $m_t' =\frac{x_1-x_0}{t_1-t_0}$ and  $\tau_t' = \eta^2\frac{t_0+t_1- 2t}{(t_1-t_0)}$.
Computing the time derivative of $p_t$ yields
\begin{align}
\partial_t p_t(x) 
&= p_t(x) \Big(
- \frac{d}{2} \frac{\tau_t'}{\tau_t} +  \frac{(x-m_t)\cdot m_t'}{\tau_t}  +\frac{\|x-m_t\|^2 \tau_t'}{2\tau_t^2}\Big)
\\
&= p_t(x) \Big(- \frac{d \, \eta^2(t_0+t_1-2t)}{ 2 \tau_t(t_1-t_0)}
+ \frac{ (x-m_t)\cdot (x_1-x_0)}{ \tau_t(t_1-t_0)} + \frac{ \|x-m_t\|^2 \eta^2(t_0+t_1-2t)}{2 \tau_t^2(t_1-t_0)} \Big)\\
&= p_t(x) \, \frac{\eta^2(t_0+t_1-2t)}{ 2 \tau_t(t_1-t_0)} \Big(-d
+ \frac{ 2(x-m_t)\cdot (x_1-x_0)}{ \eta^2 (t_0+t_1 -2t)} + \frac{ \|x-m_t\|^2 }{\tau_t} \Big) = p_t(x) \xi_t(x).\label{pt:timederiv}
\end{align}
First, we deal with the diffusion part.
Using that  $\nabla \cdot u_t = \frac{d\, \eta^2}{2\tau_t}\left(\frac{t_0+t_1-2t}{(t_1-t_0)} - 1\right)$ and
$\nabla p_t = - p_t\, \frac{x-m_t}{\tau_t}$, we can compute
\begin{multline*}
   \nabla \cdot (p_t u_t) (x) = u_t^\tT \nabla p_t + p_t \nabla \cdot u_t
   = -p_t \Big( \frac{u_t^\tT(x-m_t)}{\tau_t}  - \frac{d\, \eta^2}{2\tau_t}\left(\frac{t_0+t_1-2t}{(t_1-t_0)}-1\right)
       \Big)\\
       = -p_t \Big( \frac{(x_1-x_0)^\tT  (x-m_t)}{\tau_t(t_1-t_0)} + 
       \| x-m_t\|^2 \frac{\eta^2 }{2\tau_t^2}\left(\frac{t_0+t_1-2t}{(t_1-t_0)}-1\right) - \frac{d\, \eta^2}{2\tau_t}\left(\frac{t_0+t_1-2t}{(t_1-t_0)}-1\right)\Big).
\end{multline*}
Further, we obtain
\begin{align*}
    \Delta p_t(x) 
    &= -\nabla p_t^\tT \frac{x-m_t}{\tau_t} - \frac{d \, p_t}{\tau_t}
    = 
    p_t \Big(
    \frac{\|x-m_t\|^2}{\tau_t^2} - \frac{d}{\tau_t} \Big)
\end{align*}
so that
\begin{align}
- \nabla \cdot (p_t u_t)  + \frac12 \eta^2 \Delta p_t 
&=
p_t  \Big(\frac{
       (x_1-x_0)^\tT  (x-m_t)}{\tau_t(t_1-t_0)} + 
       \| x-m_t\|^2 \frac{\eta^2 }{2\tau_t^2}\left(\frac{t_0+t_1-2t}{(t_1-t_0)}-1\right)\\
       &\qquad- \frac{d\, \eta^2}{2\tau_t}\left(\frac{t_0+t_1-2t}{(t_1-t_0)}-1\right)       +    \frac{\eta^2}{2}\left(\frac{\|x-m_t\|^2}{\tau_t^2} - \frac{d}{\tau_t}\right) \Big)\\
       &= \Big( \frac{\eta^2}{2\tau_t}\left(-d\left(\frac{t_0+t_1-2t}{(t_1-t_0)}-1\right)-d\right) +\frac{(x_1-x_0)^\tT(x-m_t)}{\tau_t(t_1-t_0)}\\
       &\qquad+\|x-m_t\|^2 \frac{\eta^2}{2\tau_t^2}\left(\frac{t_0+t_1-2t}{(t_1-t_0)} -1 +1\right)\Big) \\
       &= \partial_t p_t(x) .
\end{align}
For the jump part, we get
\begin{align}\label{continuity:jump}
&\int p_t(y) q_t (x,y) - p_t(x) q_t(y,x){\rm d}y
= J_t(x) \int \lambda_t(y) p_t(y){\rm d}y - \lambda_t(x) p_t(x)
\end{align}
Since
$\max(0,\xi_t) = \xi_t + \lambda_t$,
we conclude by definition of $\xi_t$ that
\begin{align}
\int \max(0,\xi_t) p_t \d y &= \int \xi_t p_t \d y + \int \lambda_t p_t \d y 
= \mathbb E_{y \sim \mathcal N(m_t,\tau_t I_d)} [\xi_t] + \int \lambda_t p_t \d y\\
&= \frac{\eta^2(t_0+t_1-2t)}{2(t_1-t_0)\tau_t} (d+0-d) + \int \lambda_t p_t \d y = \int \lambda_t p_t \d y.
\end{align}
Hence, unless $\lambda_t=0$ everywhere, which is the case if and only if $t=\tfrac{t_0+t_1}{2}$ and $x_0=x_1$,
\begin{equation}\label{jump:measure}
J_t(x) = \frac{\big(\xi_t(x) + \lambda_t(x) \big) p_t(x)}{\int \lambda_t(y) p_t(y) \d y}.
\end{equation}
Inserting \eqref{jump:measure} in \eqref{continuity:jump}, and using $\lambda_t=\xi_t=0$ if $t=\tfrac{t_0+t_1}{2}$ and $x_0=x_1$, we  obtain 
\begin{align}
\int p_t(y) q_t(x,y) - p_t(x) q_t (y,x) \d y
=  p_t(x) \xi_t(x) = \partial_t p_t(x).
\end{align}
This finishes the proof.
\end{proof}

\begin{proof}[Proof of Remark \ref{prop:post_x_0}.]
    Since $Y_t^{x_0}\sim P_t^{x_0}$, it suffices to compute $P_0^{x_0}$ and $P_1^{x_0}$ as follows:
    \begin{align*}
        \int_{\R^d} f(x)\,\d P_0^{x_0}(x)&=\int_{\R^d\times \R^d} f(x)\,\d P_0(\cdot ,x_0,x_1)(x)\,\d P_{X_1|X_0=x_0}(x_1)\\
        &= \int_{\R^d \times \R^d} f(x) \,\d \NN(x_0,\rho^2 I_d)(x) \,\d P_{X_1|X_0=x_0}(x_1)=\int_{\R^d}f(x)\,\d \NN(x_0,\rho^2 I_d)(x),\\
        \int_{\R^d}f(x)\,\d P_1^{x_0}(x)&=\int_{\R^d\times \R^d}f(x)\,\d \NN(x_1,\rho^2I_d)(x)\,\d P_{X_1|X_0=x_0}(x_1)\\
        &=\int_{\R^d}f(x)\left(\int_{\R^d} ({2\pi}\rho^2)^{-d/2}e^{-\frac{\|x-x_1\|^2}{2\rho^2}}\d P_{X_1|X_0=x_0}(x_1)\right)\d x\\
        &=\int_{\R^d} f(x)\,\d \bigl(P_{X_1|X_0=x_0}*\NN(0,\rho^2 I_d)\bigr)(x). \qedhere
    \end{align*}
\end{proof}

\begin{proof}[Proof of Theorem \ref{prop2}.]
By definition of the KL-divergence, we get
\begin{align*}
    &F_{t,x}(\lambda,\mu,\Sigma)
    =  \lambda - \lambda_t(x)  +\int q_t(y,x) \log\left(\frac{q_t(y,x)}{\lambda \mathcal N(\mu, \Sigma)(y)}\right) \d y\\
    &= \lambda - \lambda_t(x)
    + \int q_t(y,x) \left(\log(q_t(y,x))-\log(\lambda) \right) \d y\\
    &\quad +  \int q_t(y,x) \left(\frac{d}{2} \log(2\pi)  + \frac{1}{2}\log(\det\Sigma) + \frac{1}{2}\left\langle y-\mu, \Sigma^{-1} (y-\mu) \right\rangle\right)\d y \\
    &= \lambda + \lambda_t(x) 
    \left( \frac{d}{2} \log(2\pi) + \frac{1}{2}\log(\det\Sigma) - \log(\lambda) 
    + \frac{1}{2} \int J_t(y) \left\langle y-\mu, \Sigma^{-1} (y-\mu) \right\rangle \d y\right)\\
    &\quad+ \int q_t(y,x)\log(q_t(y,x)) \d y - \lambda_t(x)\\
    &= \lambda + \lambda_t(x) \left( \frac{1}{2}\log(\det\Sigma) - \log(\lambda) 
    + \frac12\left\langle\Sigma^{-1},\Sigma^J\right\rangle + \frac12\left\langle \mu^J-\mu,\Sigma^{-1}(\mu^J-\mu)\right\rangle\right) + {\rm const}.
\end{align*}
Here, the identity $\int J_t(y) \left\langle y-\mu, \Sigma^{-1} (y-\mu) \right\rangle \d y = \left\langle\Sigma^{-1},\Sigma^J\right\rangle +   \left\langle \mu^J-\mu,\Sigma^{-1}(\mu^J-\mu)\right\rangle$ follows using the decomposition $y-\mu=(y-\mu^J) +(\mu^J-\mu)$.

Assume that $\lambda_t(x)\neq 0$ and $\Sigma^J$ is positive definite. Then, $\lambda_t(x)>0$ and $F_{t,x}(\lambda, \mu,\Sigma)=F_1(\lambda)+\tfrac{1}{2}\lambda_t(x)(F_2(\Sigma)+F_3(\mu,\Sigma))$, where $F_1(\lambda)=\lambda - \lambda_t(x)\log \lambda$, $F_2(\Sigma)=\log\det\Sigma + \tr(\Sigma^{-1}\Sigma^{J})$, $F_3(\mu,\Sigma)= \langle \mu^J-\mu ,\Sigma^{-1}(\mu^J-\mu)\rangle.$ Since $F_1'(\lambda)=1-\tfrac{\lambda_t(x)}{\lambda}$ and $F_1''(\lambda)=\tfrac{\lambda_t(x)}{\lambda^2}>0$, $F_1$ has a unique minimum at $\lambda=\lambda_t(x).$ Since $\Sigma^{-1}$ is positive definite, $F_3(\cdot,\Sigma)$ has a unique minimum at $\mu=\mu^J,$ whose value does not depend on $\Sigma$. Thus, it remains to show that $F_2$ has the unique minimizer $\Sigma=\Sigma^J.$ 

To this end, let $A=(\Sigma^J)^{1/2}\Sigma^{-1}(\Sigma^J)^{1/2}$. Then
\begin{equation}
\label{eq:F2}
F_2(\Sigma)-F_2(\Sigma^J) = \operatorname{tr}(A)-\log\det(A)-d = \sum_{i=1}^d \bigl(a_i-\log(a_i)-1\bigr)
\geq 0,
\end{equation}
where $a_1,\ldots,a_d>0$ are the eigenvalues of $A$. The last inequality follows since $a-\log a -1 \geq 0$ for $a>0$ with equality if and only if $a=1$. Hence, equality in \eqref{eq:F2} holds if and only if $A=I$, that is $\Sigma =\Sigma^J.$
\end{proof}

\subsection{Proof of Theorem \ref{prop1}.}\label{paragraph:prop1}
We derive closed-form expressions for the mean and covariance of $J_t$ separately for $d=1$ and $d=2$. Note that we consider arbitrary $t_0$ and $t_1$, i.e., $m_t, \tau_t, u_t, \xi_t,$ and $a_t$ are defined by the expressions at the beginning of this section.
\begin{theorem}
\label{prop:formula1d}
In dimension $d=1$, let $\mu^{ J}=\mu_{x_0,x_1,J}$ and $(\sigma^{ J})^2=(\sigma^{x_0,x_1,J})^2$ denote the mean and variance of $J_t$. 
Set \[ z_\pm \coloneqq - a_t \pm \sqrt{a_t^2 + 1}.\]
Moreover, define
\begin{align}
I_0 &\coloneqq \sqrt{\pi/2} \Big(\erf\Big(\frac{z_+}{\sqrt{2}} \Big) - \erf\Big(\frac{z_-}{\sqrt{2}} \Big)\Big), \quad
I_1 \coloneqq 
 \e^{-z_-^2/2} - \e^{-z_+^2/2}, \label{eq:I_01}\\
I_k &\coloneqq(k-1) I_{k-2} -\left(z_+^{k-1}\e^{-z_+^2/2} - z_-^{k-1} e^{-z_-^2/2}\right), \quad k=2,3,4.  \label{eq:I_234}
\end{align}
Then, it holds for $t \in \left(\frac{t_0+t_1}{2},t_1\right]$ that
\begin{align}
\mu^{J}&=  m_t + \sqrt{\tau_t} \, \frac{I_3  + 2 a_t I_2 - I_1}{I_2  + 2 a_t I_1 - I_0},
\label{eq:mean_1d} \\
(\sigma^{J})^2 &= \tau_t \frac{I_4 + 2a_t I_3 - I_2}{I_2  + 2 a_t I_1 - I_0}- (m_t - \mu^{J})^2. \label{eq:var_1d}
\end{align}
For $t \in \left[t_0,\frac{t_0+t_1}{2}\right)$, the formulas \eqref{eq:mean_1d} and \eqref{eq:var_1d} remain valid if each $I_k$ is replaced by $\tilde{I}_k \coloneqq \frac{(-1)^k+1}{2}\sqrt{2\pi}(k-1)!!-I_k$. (For $t=\frac{t_0+t_1}{2}$, $\mu^J=0$ and $\sigma^J=1$ by definition of $J_t$.)
\end{theorem}

\begin{proof}
Setting $z\coloneqq (y-m_t)/\sqrt{\tau_t}$, we obtain  
 $$ 
 \xi_t(z)   = 
\frac{\eta^2(t_0+t_1-2t)}{2\tau_t(t_1-t_0)}
\left( z^2 + 2a_t z   - 1\right).
 $$ 
 This quadratic polynomial has zeros
 $$ z_\pm = -a_t \pm \sqrt{a_t^2 + 1}. 
$$
The support of $J_t(z)$ is 
$(-\infty, z_-] \cup  [z_+,\infty)$ if $t \in [t_0,\frac{t_0+t_1}{2})$, and 
$[z_-,z_+]$ if $t \in (\frac{t_0+t_1}{2},t_1]$. For $t \in (\frac{t_0+t_1}{2},t_1]$ we set
$$
C_t \coloneqq 
\int_{z_-}^{z_+}
(z^2 + 2 a_tz -1) \,   e^{-z^2/2} \d z,
$$
and compute that
\begin{align} 
\mu^{\rm J}  &= \int y J_t(y) \d y
=  
\frac{1}{C_t}
\int_{z_-}^{z_+}
(z \sqrt{\tau_t}  + m_t) \, (z^2 + 2 a_t z -1) \,   e^{-z^2/2} \d z
\notag \\
&=
m_t + \frac{\sqrt{\tau_t}}{C_t}
\int_{z_-}^{z_+} z \,(z^2 + 2 a_t z -1) \,  \e^{-z^2/2} \d z
\label{mean}
\end{align}
and
\begin{align} \label{variance}
(\sigma^{\rm J})^2 &= \int (y- \mu^{\rm J})^2 J_t(y) \d y
=
\int y^2 J_t(y) \d y - (\mu^{\rm J})^2 \notag \\
&=
\frac{1}{C_t}  \int_{z_-}^{z_+}
(\sqrt{\tau_t} z + m_t)^2 (z^2 + 2 a_t z - 1) \, \e^{-z^2/2} \d z - (\mu^{\rm J})^2.
\end{align}
Defining
$$
I_k \coloneqq \int_{z_-}^{z_+} z^k e^{-z^2/2} \d z, \quad k=0,\ldots,4, 
$$
it follows that
\begin{align*}
\mu^{J}&=  m_t + \frac{\sqrt{\tau_t}}{C_t}(I_3  + 2 a_t I_2 - I_1)=m_t + \sqrt{\tau_t} \, \frac{I_3  + 2 a_t I_2 - I_1}{I_2  + 2 a_t I_1 - I_0},
\\
(\sigma^{J})^2
& = \frac{1}{C_t} \Big( 
\tau_t \int_{z_-}^{z_+} z^2 (z^2+ 2 a_t z -1) \,e^{-z^2/2}\,\d z 
+ 2 \sqrt{\tau_t} m_t \int_{z_-}^{z_+} z (z^2+ 2 a_t z -1)\, e^{-z^2/2}\,\d z \Big)
+ m_t^2 -(\mu^{\rm J})^2
\\
&=  \frac{\tau_t}{C_t} 
\int_{z_-}^{z_+} z^2 (z^2+ 2 a_t z -1)\, e^{-z^2/2}\,\d z  + 2 m_t (\mu^{\rm J} -m_t) + m_t^2 - (\mu^{\rm J})^2\\
&= \tau_t \frac{I_4 + 2a_t I_3 - I_2}{I_2  + 2 a_t I_1 - I_0}- (m_t - \mu^{\rm J})^2.
\end{align*}
For $t \in [t_0,\frac{t_0+t_1}{2})$, the computations are analogous using $\tilde{C}_t \coloneqq 
\int_{(-\infty, z_-] \cup  [z_+,\infty)}
(z^2 + 2 a_tz -1) \,   e^{-z^2/2} \d z$ and
$$\tilde{I}_k \coloneqq \int_{(-\infty, z_-] \cup  [z_+,\infty)} z^k e^{-z^2/2} \d z, \quad k=0,\ldots,4.$$

It remains to derive the explicit expressions for $I_k$ and $\tilde{I}_k$.
The value of $I_0$ follows by
the definition of the Gaussian error function
$$
\int_{-\infty}^z \tfrac{1}{\sqrt{2\pi}} e^{-z^2/2}\, \d z = \frac{1}{2}\left(1+{\rm erf}(z/\sqrt{2})\right).
$$
For $k \ge 1$, and with a slight abuse of notation for $k=1$, integration by parts gives
$$
\int_{z_-}^{z_+} z^k \e^{-z^2/2} \d z = - z^{k-1} \e^{-z^2/2}|^{z_+}_{z_-} +(k-1)\int_{z_-}^{z_+} z^{k-2} \e^{-z^2/2} \d z,
$$
so that
\begin{align*}
I_k &=(k-1) I_{k-2} -\left(z_+^{k-1}\e^{-z_+^2/2} - z_-^{k-1} e^{-z_-^2/2}\right).
\end{align*}
This yields \eqref{eq:I_01} and \eqref{eq:I_234}.
Using the well-known formula for Gaussian moments 
\begin{equation}
   \int_{-\infty}^\infty z^k e^{-z^2/2} {\rm d}z/\sqrt{2\pi} = \begin{cases} 0 & \quad k\mbox{ odd,}\\ (k-1)!!&\quad k\mbox{ even,}\end{cases}
\end{equation}
we finally obtain that $\tilde{I}_k =\frac{(-1)^k+1}{2}\sqrt{2\pi}(k-1)!!-I_k$. 
\end{proof}

\begin{theorem}\label{prop:2d_formulas}
Let $d=2$ and $t\neq \tfrac{t_0 +t_1}{2}$. Set
\begin{equation*}
    f_k(r,a_t)\vcentcolon= \int_0^{2\pi} \cos(\phi)^k \max(0,r^2+2\Vert a_t\Vert r\cos \phi -d) \,\textup{d}\phi\qquad \textup{for }r\geq 0 \textup{ and }k=0,1,2.
\end{equation*}
Define
\begin{align*}
    D &= \int_0^\infty e^{-\frac{r^2}{2}}r f_0(r, a_t){\rm d}r, &N_1&= \int_0^\infty e^{-\frac{r^2}{2}} r^2 f_1(r,a_t) {\rm d}r, \\
    N_{11} &= \int_0^\infty e^{-\frac{r^2}{2}} r^3 f_2(r,a_t) {\rm d}r, & N_{r^2} &= \int_0^\infty e^{-\frac{r^2}{2}} r^3 f_0(r,a_t) {\rm d}r.
\end{align*}
Then,
\begin{align}\label{2dformulas}
    \mu^J = E[J_t] &=  m_t + \sqrt{\tau_t}\frac{N_1}{D}\frac{a_t}{\|a_t\|}, \notag\\ 
    \Cov(J_t) &= \tau_t \left( \frac{2N_{11}-N_{r^2}}{D}-  \frac{N_1^2}{D^2}\right)\frac{a_t a_t^\tT }{\|a_t\|^2} 
    + \tau_t \frac{N_{r^2}-N_{11}}{D} I_d.
\end{align}
\end{theorem}
Note that the $f_k$ above can be computed in closed form. Thus, even though our notation suggests otherwise, the computation of $\mu^J$ and $\Cov(J_t)$ does not require the evaluation of iterated integrals.

\begin{proof}
Let $R_t\in\R^{d\times d}$ be an orthogonal matrix with ${\rm det}(R_t)=1$ and $R_t a_t = \|a_t\| e_1$, where $e_1\in\R^d$ denotes the first cartesian unit vector. With the substitution $v \vcentcolon= S_t(y) \vcentcolon= R_t\frac{y-m_t}{\sqrt{\tau_t}}$, we obtain that
\begin{equation*}
    \quad \xi_t(y)=\kappa_t \left(\frac{\|y-m_t\|^2}{\tau_t}+2 a_t^\tT  \frac{y-m_t}{\sqrt{\tau_t}} - d\right) 
    = \kappa_t (\Vert v\Vert^2+2 \Vert a_t\Vert v_1 -d) \quad\textup{with}  \quad \kappa_t \vcentcolon= \frac{\eta^2(t_0 +t_1 -2t)}{2\tau_t(t_1-t_0)},
\end{equation*}
and the density $p_t$ of $\mathcal{N}(m_t,\tau_t I_d)$ takes the form $p_t(y)=Z_t e^{-\Vert v\Vert^2/2}$ for some positive constant $Z_t$. Consequently,
\begin{align*}
\max(0,\xi_t(y))p_t(y){\rm d}y &= \kappa_t Z_t\cdot \max\left(0, \|v\|^2 + 2 \|a_t\| v_1 - d \right) \cdot e^{-\|v\|^2/2} {\tau_t}^{d/2} {\rm d}v.
\end{align*}

Define the moments of $J_t$ (up to a constant factor) as
\begin{equation*}
    M_k \vcentcolon= \int_{\R^2} y^{\otimes k}\max(0,\xi_t(y))p_t(y),\d y,
\qquad k=0,1,2,
\end{equation*}
where $y^{\otimes 0}=1$, $y^{\otimes 1}=y$, and $y^{\otimes 2}=yy^\tT$.
Transforming to polar coordinates we compute that 
\begin{align*}
    M_0 
    &= \kappa_t Z_t \tau_t \int_{\R^2} \max\left(0, \|v\|^2 + 2 \|a_t\| v_1 - d \right) \, e^{-\|v\|^2/2} {\rm d}v \\
    &= \kappa_t Z_t \tau_t \int_0^\infty e^{-r^2/2} r \int_0^{2\pi} \max\left(0, r^2 + 2 \|a_t\| r\cos\phi - d \right) {\rm d}\phi \,\d r \\
    &= \kappa_t Z_t \tau_t \int_{0}^\infty e^{-r^2/2} r \, f_0(r,a_t) \,\d r.
\end{align*}
Since $y=\sqrt{\tau_t}R_t^\tT v+m_t$, we obtain that
\begin{align*}
    M_1 
    &= \kappa_t Z_t \tau_t^{3/2} R_t^\tT  \int_{\R^2} v \max\left(0, \|v\|^2 + 2 \|a_t\| v_1 - d \right) \, e^{-\|v\|^2/2} \d v+ m_t M_0 \\
    &= \kappa_t Z_t \tau_t^{3/2} \frac{a_t}{\Vert a_t\Vert} \int_{0}^\infty  r^2 e^{-r^2/2} f_1(r,a_t) \d r+ m_t M_0,
\end{align*}
where we used that by symmetry
\begin{align*}
   \int_{\R^2} v &\max\left(0, \|v\|^2 + 2 \|a_t\| v_1 - d \right) \, e^{-\|v\|^2/2} \d v \\
   &= e_1 \int_{\R^2} v_1 \max\left(0, \|v\|^2 + 2 \|a_t\| v_1 - d \right) \, e^{-\|v\|^2/2} \d v \\
   &= e_1 \int_0^\infty r^2 e^{-r^2/2}\int_0^{2\pi} \cos{\phi} \max(0,r^2 + 2 \Vert a_t \Vert r \cos\phi - d) \,\d \phi \, \d r.
\end{align*}
Similarly, for $M_2$ we compute with $y=\sqrt{\tau_t}R_t^\tT v+m_t$ that
\begin{align*}
    M_2 &= \kappa_t Z_t \tau_t^{2} R_t^\tT  \int_{\R^2} vv^\tT  \max\left(0, \|v\|^2 + 2 \|a_t\| v_1 - d \right) \, e^{-\|v\|^2/2} \d v \, R_t \\
    &\qquad+ M_1 m_t^\tT  + m_t M_1^\tT - m_tm_t^\tT  M_0.
\end{align*}
Here, by symmetry,
\begin{align*}
    \int_{\R^2} vv^\tT  &\max\left(0, \|v\|^2 + 2 \|a_t\| v_1 - d \right) \, e^{-\|v\|^2/2} \d v \\
    &= \int_{\R^2} (v_1^2 (e_1e_1^\tT ) + v_2^2(e_2 e_2^\tT ))  \max\left(0, \|v\|^2 + 2 \|a_t\| v_1 - d \right) \, e^{-\|v\|^2/2} \d v  \\
    &= \int_{\R^2} ((v_1^2 -v_2^2) (e_1e_1^\tT ) + v_2^2 I_d)  \max\left(0, \|v\|^2 + 2 \|a_t\| v_1 - d \right) \, e^{-\|v\|^2/2} \d v.
\end{align*}
Hence, 
\begin{align*}
    R_t^\tT  \int_{\R^2} vv^\tT  &\max\left(0, \|v\|^2 + 2 \|a_t\| v_1 - d \right) \, e^{-\|v\|^2/2} \d v \, R_t \\
    &= \frac{a_ta_t^\tT }{\Vert a_t\Vert^2} \int_{\R^2} (v_1^2-v_2^2)\max\left(0, \|v\|^2 + 2 \|a_t\| v_1 - d \right) \, e^{-\|v\|^2/2} \d v \\
    &\quad + I_d \int_{\R^2} v_2^2\max\left(0, \|v\|^2 + 2 \|a_t\| v_1 - d \right) \, e^{-\|v\|^2/2} \d v, \\
    &= \frac{a_ta_t^\tT }{\Vert a_t\Vert^2} \int_{0}^\infty r^3 e^{-r^2/2} \int_0^{2\pi}  (2\cos(\phi)^ 2 -1) \max\left(0, r^2 + 2 \|a_t\| r \cos \phi - d \right) \d \phi \,\d r \\
    &\quad + I_d \int_{0}^\infty  r^3 e^{-r^2/2} \int_0^{2\pi}  (1-\cos(\phi)^2) \max\left(0, r^2 + 2 \|a_t\| r \cos \phi - d \right) \d \phi \, \d r,\\
    &= \frac{a_ta_t^\tT }{\Vert a_t\Vert^2} \int_{0}^\infty r^3 e^{-r^2/2} (2f_2(r, a_t) - f_0(r, a_t)) \d r \\
    &\quad + I_d \int_{0}^\infty r^3 e^{-r^2/2} (f_0(r,a_t) - f_2(r, a_t)) \d r.
\end{align*}
Finally, the claim follows since $E[J_t]=M_1/M_0$ and $\Cov(J_t)=M_2/M_0 - E[J_t]E[J_t]^\tT $.
\end{proof}

\section{Proof of Theorem \ref{rem:joint_post}}\label{app:joint_post}
In order to derive the stability result, we proceed in two steps. First, we show that we can approximate $P_{X_0,X_1}$ if we can approximate $P_{X_0}$ and $P_{X_1|X_0=x_0}$ by $P_{Y_0}$ and a Markov kernel $\mu(\cdot,x_0)$, respectively.

\begin{lemma}\label{lemma:disintegration_lipschitz}	Let $P_{X_0,X_1}\in \mathcal P_2(\R^m\times\R^d)$, $P_{Y_0}\in\PP_2(\R^m)$ and let $\mu:\mathcal B(\R^d)\times \R^m\to [0,1]$ be a Markov kernel. Assume that
    \begin{romanenumerate}
	    \item $\int W_2^2(P_{X_1|X_0=a},\mu(\cdot,a))\,\d P_{X_0}(a)\leq \epsilon$,
    \item $W_2(\mu(\cdot,a),\mu(\cdot,b))\leq K\|a-b\|$.
    \end{romanenumerate}
    Then $\alpha\coloneqq \mu(\cdot,x_0)\times_{x_0}P_{Y_0}\in \PP_2(\R^m\times\R^d)$ satisfies
    \begin{align}
    W_2^2(P_{X_0,X_1},\alpha)\leq (1+2K^2)W_2^2(P_{X_0},P_{Y_0}) + 2 \varepsilon.
    \end{align}
\end{lemma}
\begin{proof}
    Let $\gamma_{a,b}\in \Gamma_o(P_{X_1|X_0=a},\mu(\cdot,b))$ and let $\beta\in \Gamma_o(P_{X_0},P_{Y_0})$. Note that $\gamma_{a,b}\times_{a,b} \beta\in\Gamma(P_{X_0,X_1},\alpha)$. Therefore,
\begin{align*}
&W_2^2(P_{X_0,X_1},\alpha)\leq \int \|(a,x)-(b,y)\|^2\d\gamma_{a,b}(x,y)\,\d \beta(a,b)\\
&~=\int \|a-b\|^2\d\beta + \int \|x-y\|^2\d\gamma_{a,b}(x,y)\,\d\beta(a,b)\\
&~=W_2^2(P_{X_0},P_{Y_0}) + \int W_2^2\left(P_{X_1|X_0=a},\mu(\cdot,b)\right)\d\beta(a,b) \\
&~\leq W_2^2(P_{X_0},P_{Y_0})+ \int \left(W_2(P_{X_1|X_0=a},\mu(\cdot,a))+W_2(\mu(\cdot,a),\mu(\cdot,b))\right)^2\d\beta(a,b)\\
&~\leq W_2^2(P_{X_0},P_{Y_0})+ 2 \int W_2^2(P_{X_1|X_0=a},\mu(\cdot,a))+W^2_2(\mu(\cdot,a),\mu(\cdot,b))\,\d\beta(a,b)\\
&~\leq W_2^2(P_{X_0},P_{Y_0})+ 2 \varepsilon + 2 \int K^2 \Vert a-b\Vert^2 \,\d\beta(a,b)\\
&~= (1+2K^2)W_2^2(P_{X_0},P_{Y_0}) + 2 \varepsilon. \qedhere
\end{align*}
\end{proof}

Inductive application of the previous lemma, yields the following precise version of Theorem \ref{rem:joint_post}.

\begin{theorem}
    \label{thm:stability_exact}
    Let $X_0,\ldots,X_n, Y_0\in L^2(\Omega,\mathbb P)$ and, for $i=0,\ldots,n-1$, let $\mu_i(\cdot,x_0,\ldots,x_i)$ be Markov kernels in $x_0,\ldots,x_i$. Assume that
    \begin{romanenumerate}
        \item $\int W_2^2(P_{X_{i+1}|X_0=a_0,\ldots ,X_i=a_i},\mu_i(\cdot,a_0,\ldots,a_i))\,\d P_{X_0,\ldots,X_i}(a_0,\ldots,a_i) \leq \epsilon$,
        \item $W_2\left(\mu_i(\cdot,a_0,\ldots,a_i),\mu_i(\cdot,b_0,\ldots,b_i)\right)\leq K\|(a_0,\ldots,a_i)-(b_0,\ldots,b_i)\|$.
    \end{romanenumerate}
    Then, $\hat{P} \coloneqq \mu_{n-1}(\cdot,x_0,\ldots,x_{n-1})\times_{x_0,\ldots,x_{n-1}}\cdots\times_{x_0}P_{Y_0}\in \PP_2(\R^{(n+1)d})$ satisfies
    \begin{align}
        W_2^2(P_{X_0,\ldots,X_n},\hat{P}) \leq 2\varepsilon\sum_{i=0}^{n-1} (1+2K^2)^i+ (1+2K^2)^n \, W_2^2(P_{X_0},P_{Y_0}).
    \end{align}
\end{theorem}

Note that in Section \ref{section:whole_times_series} the conditional generators are trained to  approximate convolutions $P_{X_{i+1}|X_0=x_0,\ldots,X_i=x_i}*\NN(0,\rho^2 I_d)$, not the exact posteriors. By \cite[Lemma 7.1.10]{ambrosio2008gradient} the resulting approximation error is $$W_2(P_{X_{i+1}|X_0=x_0,\ldots,X_i=x_i},P_{X_{i+1}|X_0=x_0,\ldots,X_i=x_i}*\NN(0,\rho^2 I_d))\leq \sqrt{d}\rho.$$ Theorem \ref{thm:stability_exact} gives an estimate how this error propagates. (Note that in practice there will be an additional approximation error as the trained generators will not approximate the smoothed posterior perfectly.)

\section*{Funding}
This work was supported by the Deutsche Forschungsgemeinschaft (DFG, German Research
Foundation) under Germany´s Excellence Strategy – The Berlin Mathematics
Research Center MATH+ (EXC-2046/1, EXC-2046/2, project ID: 390685689).

\section*{Acknowledgments}
Generative AI tools were used for assistance with coding, language proofreading, and stylistic suggestions.

\section*{Data Availability}
The code used to generate the numerical results in this article is not publicly available at this time. A cleaned and documented implementation is planned for a future release.

\bibliographystyle{plain}
\bibliography{references.bib}

\beginsupplement

\begin{center}
{\Large\bfseries Supplementary Material\par}
\medskip
\end{center}

In this supplementary material, we collect additional definitions used in the paper and present further numerical results and implementation details.

\section{Additional Definitions} \label{sec:supp:defs}
Let $\mathcal M_{+}^{\text{ac}}(\R^d)$ denote the space of finite non-negative measures on $\R^d$ which are  absolutely continuous with respect to the Lebesgue measure.
Our loss function for the jump process relies  on the \emph{Kullback-Leibler  divergence} 
$\KL: \mathcal M_{+}^{\text{ac}}(\R^d) \times  \mathcal M_{+}^{\text{ac}}(\R^d) \to [0,\infty]$
defined by
\begin{align*}
    \KL(p,q):=\int  p(x) \log\Big( \frac{p(x)}{q(x)} \Big)
    - p(x)  + q(x) 
    \d x.
\end{align*}
if $p(x) = 0$ whenever $q(x) = 0$ for a.e. $x \in \R^d$, and
$\KL(p,q) \coloneqq \infty$ otherwise. We have that $\KL(p,q)=0$
if and only if $p=q$ a.e..

For a measure $\alpha \in\mathcal P(\R^d \times \R^d)$ 
and projections $\pi^x,\pi^y: \R^d \times \R^d \to \R^d$ defined by
$$\pi^x(x,y) \coloneqq x, \quad \pi^y(x,y) \coloneqq y$$
we have that
$\pi^x_\sharp \alpha$ and $\pi^y_\sharp \alpha$
are the left and right marginals of $\alpha$, respectively.
Then for $\alpha \in\mathcal P(\R^d \times \R^d)$ with marginal $\pi^x_\sharp \alpha=\mu$, there exists a 
$\mu$-a.e. uniquely defined
family of probability measures $\{\alpha_{x}\}_{x}$,  called 
\emph{disintegration of $\alpha$ with respect to  $\pi^x$}, such that 
the map $x \mapsto \alpha_{x}(B)$ is measurable for every  $B\in \mathcal B(\R^d)$, and 
$$
\alpha = \alpha_{x} \times_x \mu
$$
meaning that
\[
\int_{\R^d \times \R^d} f(x,y)  \d \alpha(x,y)
=
\int_{\R^d} \int_{\R^d} f(x,y) \d \alpha_{x}(y) \d \mu(x)
\]
for every measurable, bounded function $f:\R^d\times \R^d \to \mathbb R$.
Similarly, we define for a measure $\alpha \in\mathcal P(\R^d \times \R^d)$ with marginal $\pi^y_\sharp \alpha=\nu$
the \emph{disintegration of $\alpha$ with respect to  $\pi^y$} as
$$
\alpha = \alpha_{y} \times_y \nu.
$$
The notation of disintegration is directly related to  Markov kernels.
A \emph{Markov kernel} is a map $\mathcal K:\R ^d \times \mathcal B(\R^d)\to \R$ such that
    \begin{romanenumerate}
    \item $\mathcal K(x,\cdot)$ is a probability measure on $\R^d$ for every $x \in \R^d$, and  
    \item $\mathcal K(\cdot,B)$ is a Borel measurable map for every $B\in\mathcal B(\R^d)$.
    \end{romanenumerate}
Hence, given a probability measure $\mu \in \mathcal P(\R^d)$, 
we can define a new measure $\alpha \coloneqq \alpha_x \times_x \mu\in \mathcal P(\R^d \times \R^d)$ by 
\[
\int f(x,y)  \d \alpha (x,y) \coloneqq \int \int f(x,y) \, \d \mathcal K(x,\cdot)(y)\d \mu( x)
\]
for all measurable, bounded functions $f$.
Identifying $\alpha_{x}(B)$ with $\mathcal K(x,B)$, 
we see that conversely,
$\{\alpha_{x}\}_{x}$ is the disintegration of $\alpha$ with respect to $\pi^x$.
 Let $X_0,X_1:\Omega \to \R^d$ be random variables with joint distribution $ P_{X_0,X_1}$.
Then the conditional distribution $\{  P_{X_1|X_0=x} \}_x$ provides the disintegration 
of $P_{X_0,X_1}$, i.e.,
$$
\int f(x,y) \, \d P_{X_0,X_1}(x,y) 
= \int \int  f(x,y) \, \d P_{X_1|X_0=x} (y) \d P_{X_0}(x). 
$$
In other words, $\mathcal K = P_{X_1|X_0}$ is a Markov kernel.
If  $P_{X_0,X_1}$ admits a density $p_{X_0,X_1}$ and $P_{X_0}$ a density $p_{X_0}>0$, then $P_{X_1|X_0=x}$ has the density 
$p_{X_1|X_0 = x} = p_{X_0,X_1}(x, \cdot)/p_{X_0}(x)$.

\section{Additional Details on the Experiment with Black-Scholes Data}
\label{sec:supp:stock_details}

\paragraph{Implementation Details}
We train each method in each observation regime using 5 different training seeds. We fix $\rho^2 =0.001$ and memory length $m=20$. The hyperparameters $\eta^2$ and, for the Markov superposition, $\alpha$ are tuned on the grids $\{0.1, 0.3, 1, 3, 10\}$ and $\{0.1, 0.2\ldots, 0.9\}$. 

The grid for $\eta^2$ is chosen to cover the relevant scale of the Black-Scholes noise: The underlying geometric Brownian motion from which the one-dimensional Black-Scholes data is sampled is given by
\begin{equation*}
    \textup{d}X_t = 2 X_t \, \textup{d}t + 0.3 X_t \,\textup{d}W_t, \qquad X_0=1.
\end{equation*}
Therefore, the local diffusion coefficient $0.3x$ ranges from $0.3$ at $x=1$ to $3$ at $x=10$. Our grid for $\eta^2$ covers this range as it corresponds to $\eta \in \{0.316, \ldots, 3.162\}$.
The two-dimensional Black-Scholes data exhibits two independent components. The first corresponds to the geometric Brownian motion above, the second to the following:
\begin{equation*}
    \textup{d}Y_t = -0.2 Y_t \, \textup{d}t + 0.1 Y_t \,\textup{d}B_t, \qquad Y_0=10.
\end{equation*}
Again, the local diffusion coefficient observed for typical realizations on $[0,1]$ lies in a range covered by the grid above.

For the generator-matching models and the TFM baseline, we use standard feedforward MLPs with ReLU activations and with four hidden layers of width $256$. We train using a learning rate of $10^{-5}$ for $500$ epochs. (Note that per epoch, from each trajectory only one pair of adjacent observation points is used for training.) Autoregressive generation starts from the respective initial value $x_0$ with empty memory, i.e. memory fully padded with the repeated initial values $(x_0,0)$. The SDE matching baseline uses the architecture and generation mechanics from \cite{bartosh2025sde}. We train this baseline for $300$ epochs with learning rate $10^{-4}$.
For all models, we evaluate the MMD between generated and ground truth samples (from our validation set) as validation loss after every epoch and save the network with the best validation loss.

\paragraph{Additional Ablations} In addition to the results reported in Section \ref{sec:bs_data} of the article, we provide two ablations, investigating the influence of the chosen memory representation and the smoothing parameter $\rho^2.$ Table \ref{tab:mmd_subsampling:stock1d_zero_memory} shows model performance, when no memory is supplied, i.e., when the network inputs are only the current position and time, and the initial and terminal time of the current interval. Since Black-Scholes dynamics are Markovian, no benefit is expected from supplying a memory in this setting. Indeed, model performance is comparable in the 10\%-observation regime, and improves slightly in the 25\%- and 100\%-regime compared to Table \ref{tab:mmd_subsampling:stock}. Note that, with memory length $m=20$, in the $10\%$-observation regime at most half the memory can ever contain non-padding values, whereas the full memory can contain non-padding values in the 25\%- and 100\%-regime. The improved performance in these regimes when no memory is supplied suggests that the flat memory representation used here may not be optimal when many observations are available. This identifies one reason why model performance can decrease in our experiments as the number of observation points increases.

\begin{table}
  \centering
  \begin{minipage}{\linewidth}
  \centering
  \scalebox{0.59}{
  \begin{tabular}{llll}
  \hline
  \textbf{Method} & \textbf{10\% observed data} & \textbf{25\% observed data} & \textbf{100\% observed data} \\
        \hline
        Jump-based method
        & [0.016 {\scriptsize$\pm$ 0.0048}] {\scriptsize$(\eta^2=0.3)$}
        & [\textbf{0.0076} {\scriptsize$\pm$ 0.0021}] {\scriptsize$(\eta^2=0.3)$}
        & [\textbf{0.010} {\scriptsize$\pm$ 0.0056}] {\scriptsize$(\eta^2=0.3)$} \\

        SDE-based method
        & [0.0078 {\scriptsize$\pm$ 0.0010}] {\scriptsize$(\eta^2=1.0)$}
        & [\textbf{0.0084} {\scriptsize$\pm$ 0.0011}] {\scriptsize$(\eta^2=1.0)$}
        & [\textbf{0.0092} {\scriptsize$\pm$ 0.0018}] {\scriptsize$(\eta^2=1.0)$} \\

        Jump + SDE (MS)
        & [\textbf{0.0056} {\scriptsize$\pm$ 0.0007}] {\scriptsize$(\eta^2=1.0,\alpha=0.9)$}
        & [\textbf{0.0084} {\scriptsize$\pm$ 0.0023}] {\scriptsize$(\eta^2=0.3,\alpha=0.1)$}
        & [\textbf{0.0076} {\scriptsize$\pm$ 0.0032}] {\scriptsize$(\eta^2=0.3,\alpha=0.1)$} \\
        \hline
        \end{tabular}
        }
  \vspace{1ex}
  \caption{Test MMD over 5 runs for the no-memory ablation on the one-dimensional Black-Scholes dataset. The selected model noise parameter is shown in parentheses; for Markov
  superpositions, $\alpha$ denotes the SDE weight.}
  \label{tab:mmd_subsampling:stock1d_zero_memory}
  \end{minipage}
  \end{table}

A second reason is given by the observation that the optimal length of Markov bridge intervals is expected to depend on the noise level present in the data. Due to the diffusive scaling of the Black-Scholes dynamics, it becomes increasingly difficult to recover the local drift from subsequent observations $x_t$ and $x_{t+\Delta t}$ as $\Delta t\to 0$. 

Moreover, we train generators to approximate smoothed versions $P_{X_{i+1}|X_0=x_0,\ldots,X_i=x_i}*\mathcal{N}(0,\rho^2 I_d)$ of the underlying conditional distributions. So, with the number of bridges, the approximation errors due to the smoothing alone can accumulate, see also Theorem \ref{thm:stability_exact}. Nevertheless, in the present case, reducing $\rho^2$ does not improve performance consistently, see Table \ref{tab:mmd_stock1d_rho_ablation_101}. Moreover, for smaller $\rho^2$ the optimal hyperparameter $\eta^2$ begins to fluctuate.

\begin{table}
\centering
\begin{minipage}{\linewidth}
\centering
\scalebox{0.59}{
\begin{tabular}{lll}
\hline
\textbf{Method} & \textbf{$\rho^2=10^{-3}$} & \textbf{$\rho^2=10^{-5}$} \\
      \hline
      Jump-based method
      & [0.015 {\scriptsize$\pm$ 0.0057}] {\scriptsize$(\eta^2=0.3)$}
      & [0.019 {\scriptsize$\pm$ 0.0055}] {\scriptsize$(\eta^2=3.0)$} \\

      SDE-based method
      & [0.027 {\scriptsize$\pm$ 0.0061}] {\scriptsize$(\eta^2=0.3)$}
      & [0.019 {\scriptsize$\pm$ 0.0070}] {\scriptsize$(\eta^2=1.0)$} \\

      Jump + SDE (MS)
      & [0.012 {\scriptsize$\pm$ 0.0069}] {\scriptsize$(\eta^2=0.3,\alpha=0.3)$}
      & [0.016 {\scriptsize$\pm$ 0.0045}] {\scriptsize$(\eta^2=1.0,\alpha=0.9)$} \\

      \hline
      \end{tabular}
      }
\vspace{1ex}
\caption{Test MMD over 5 runs on the one-dimensional Black-Scholes dataset in the full observation regime for the default value and a smaller value of the smoothing parameter $\rho^2$.}
\label{tab:mmd_stock1d_rho_ablation_101}
\end{minipage}
\end{table}

\paragraph{Additional Illustrations for One-Dimensional Experiments}
Below, in addition to Figure \ref{fig:Black-Scholes-sub25} in the article, we provide more trajectory plots corresponding to the experiments reported in Table \ref{tab:mmd_subsampling:stock}.

\begin{figure}[H]
    \centering
    \includegraphics[width=0.5\linewidth]{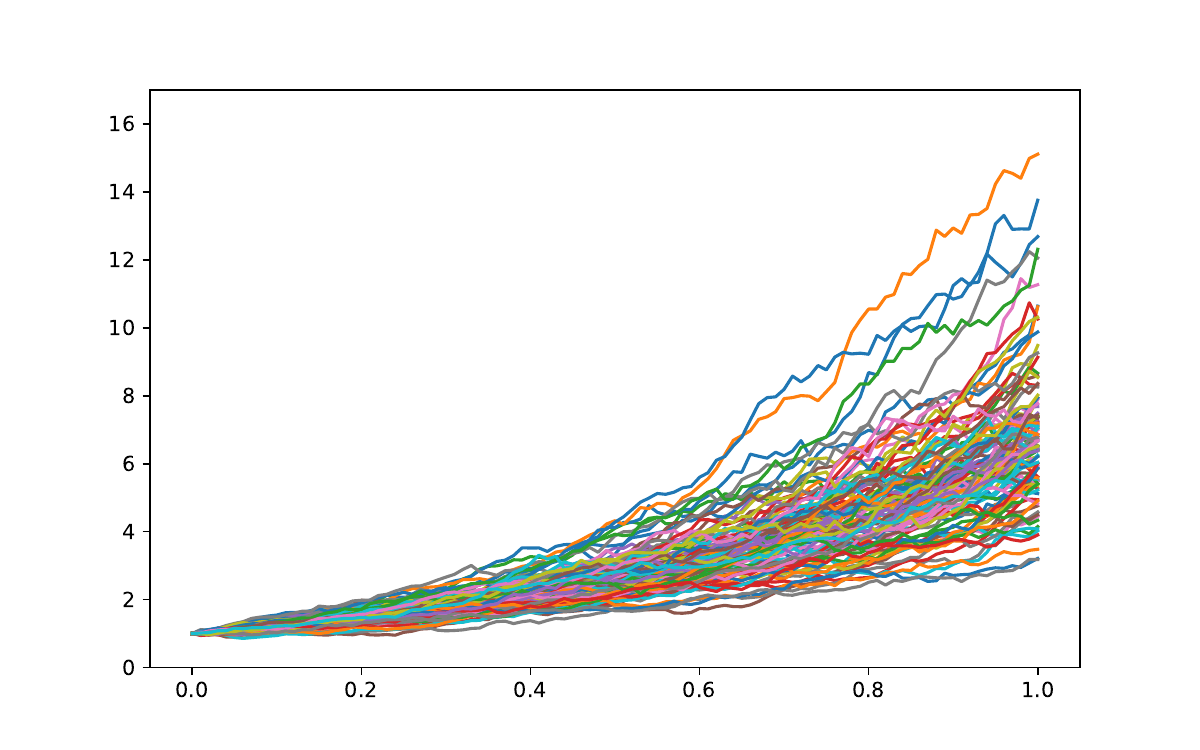}
    \caption{Ground truth trajectories of one-dimensional Black-Scholes dataset.}
\end{figure}

\begin{figure}[H]
    \centering
    \includegraphics[width = 0.32\textwidth, trim=15mm 3mm 18mm 3mm,
    clip]{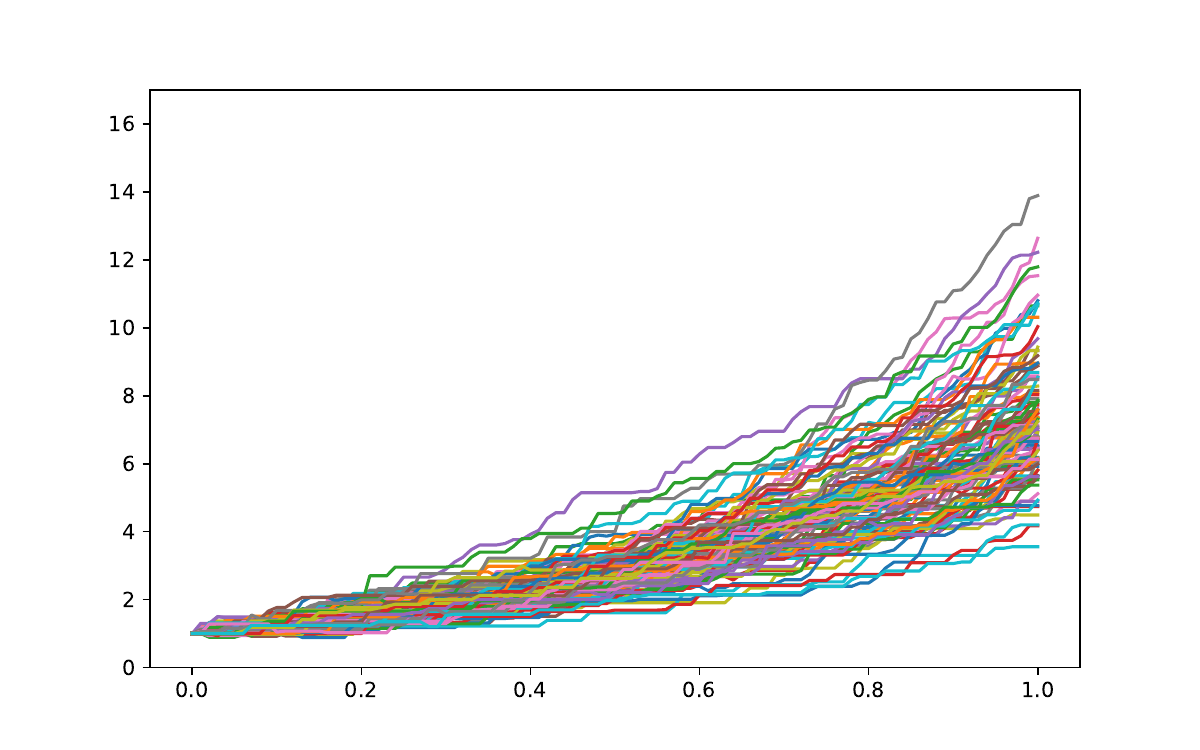}
    \includegraphics[width = 0.32\textwidth, trim=15mm 3mm 18mm 3mm,
    clip]{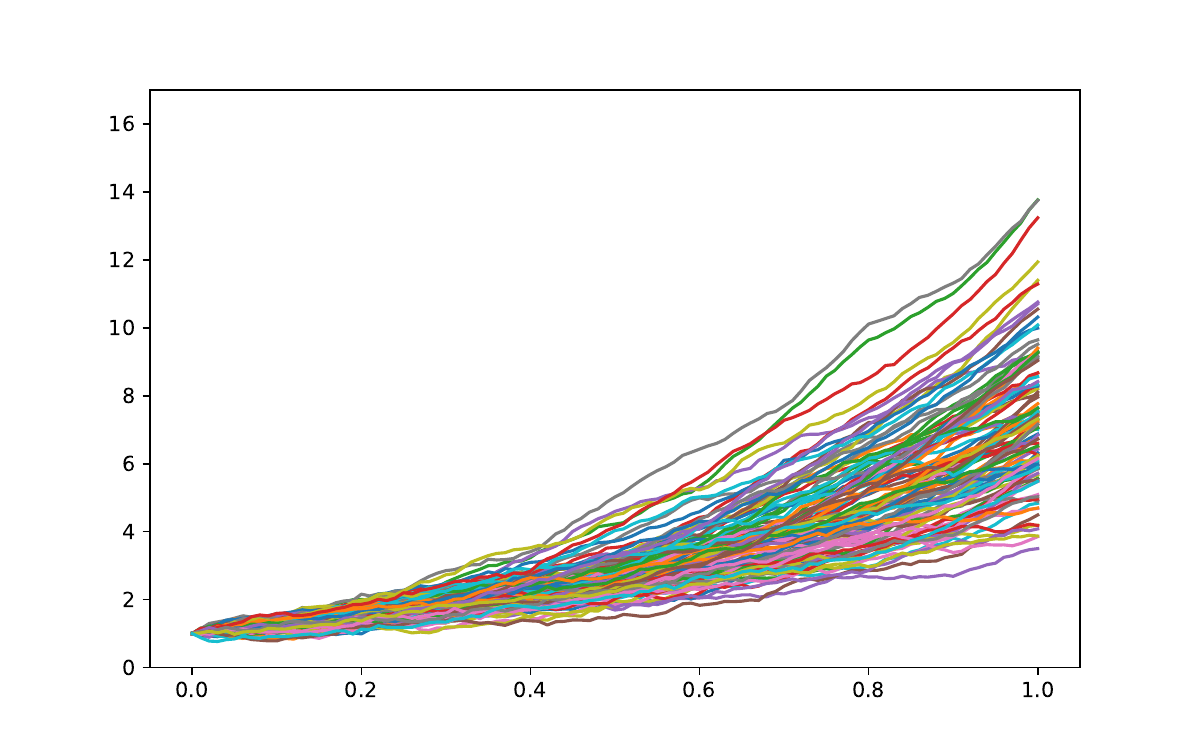}
    \par
    \includegraphics[width = 0.32\textwidth, trim=15mm 3mm 18mm 3mm,
    clip]{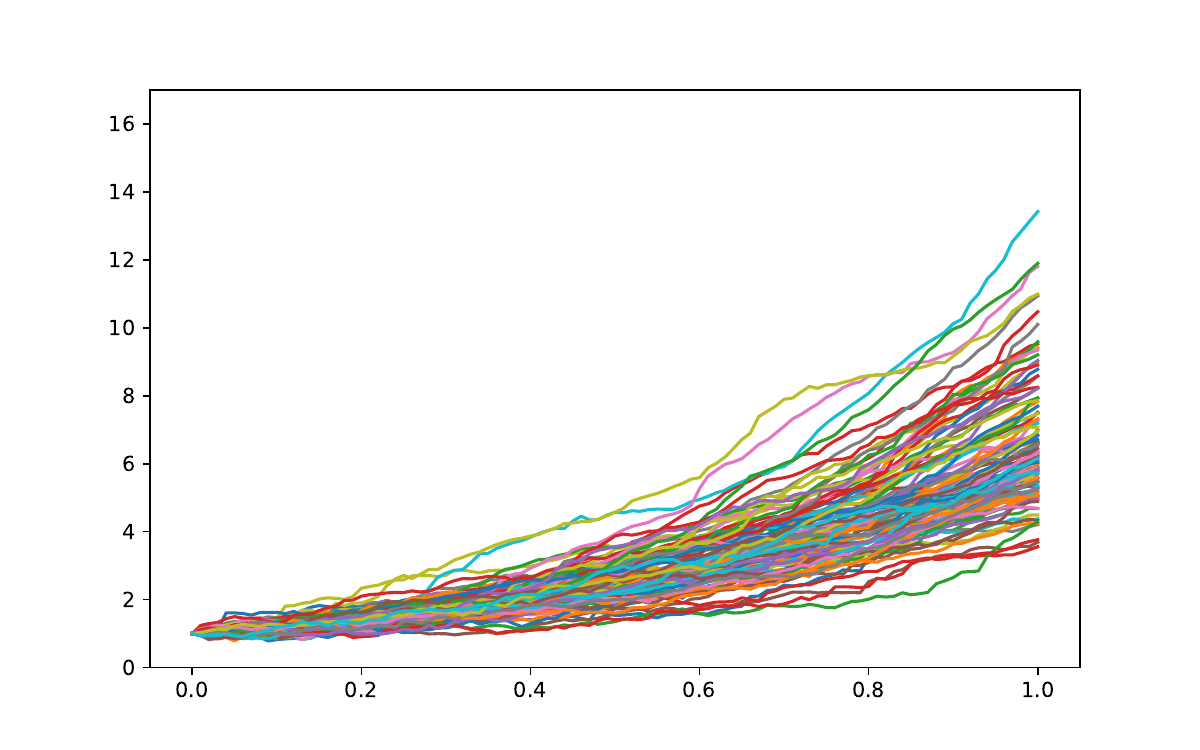}
    \includegraphics[width = 0.32\textwidth, trim=15mm 3mm 18mm 3mm,
    clip]{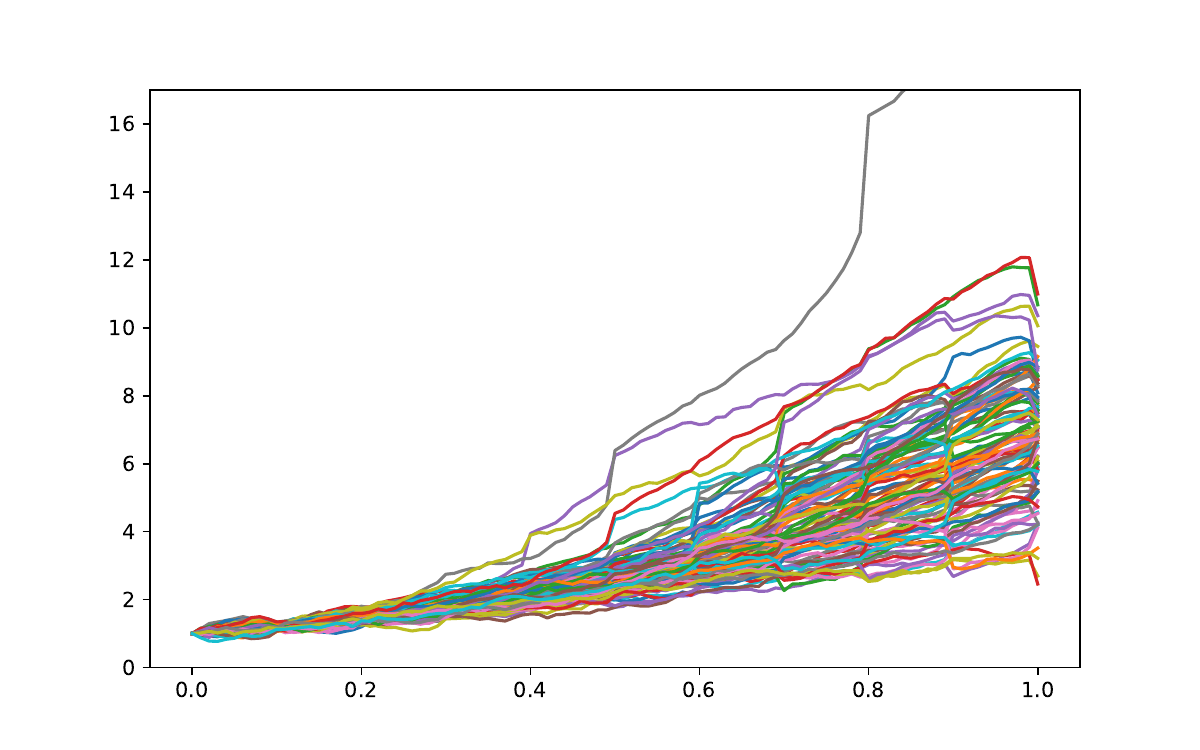}
    \includegraphics[width = 0.32\textwidth, trim=15mm 3mm 18mm 3mm,
    clip]{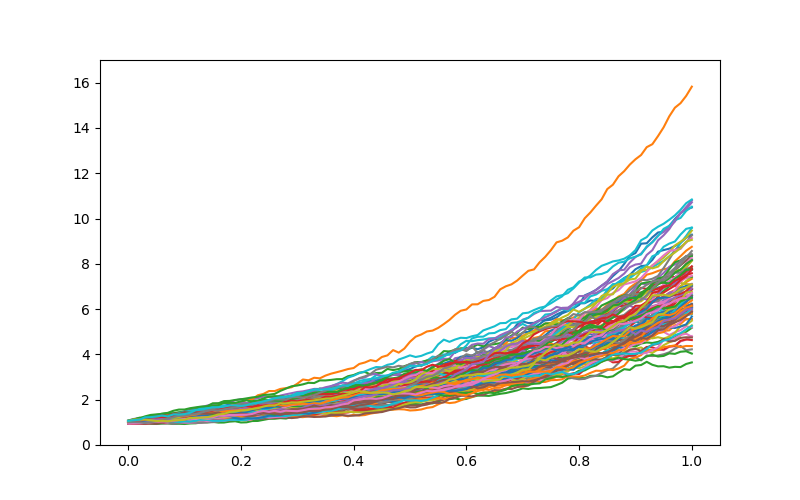}
     \caption{Trajectories generated by the following models trained in the 10\% observation regime. From left to right, top to bottom: JUMP, SDE,  MS ($\alpha=0.9$), TFM, SDE matching.}
\end{figure}

\begin{figure}[H]
    \centering
    \includegraphics[width = 0.32\textwidth, trim=15mm 3mm 18mm 3mm,
    clip]{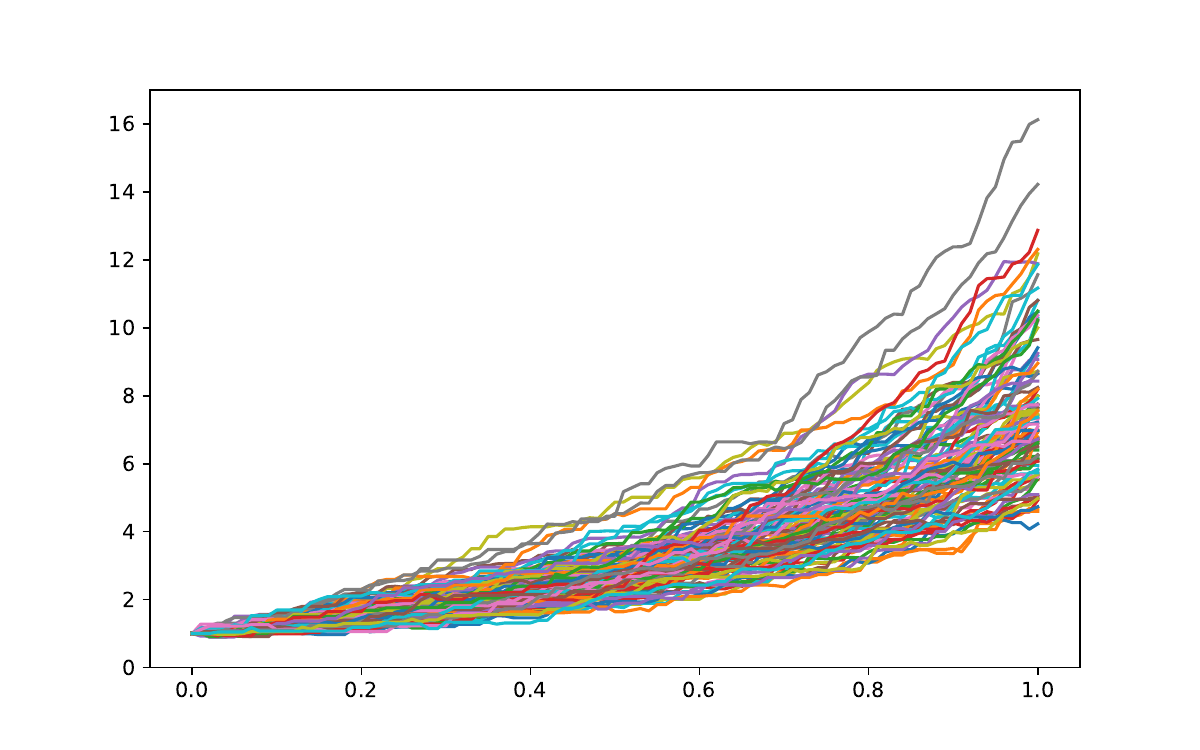}
    \includegraphics[width = 0.32\textwidth, trim=15mm 3mm 18mm 3mm,
    clip]{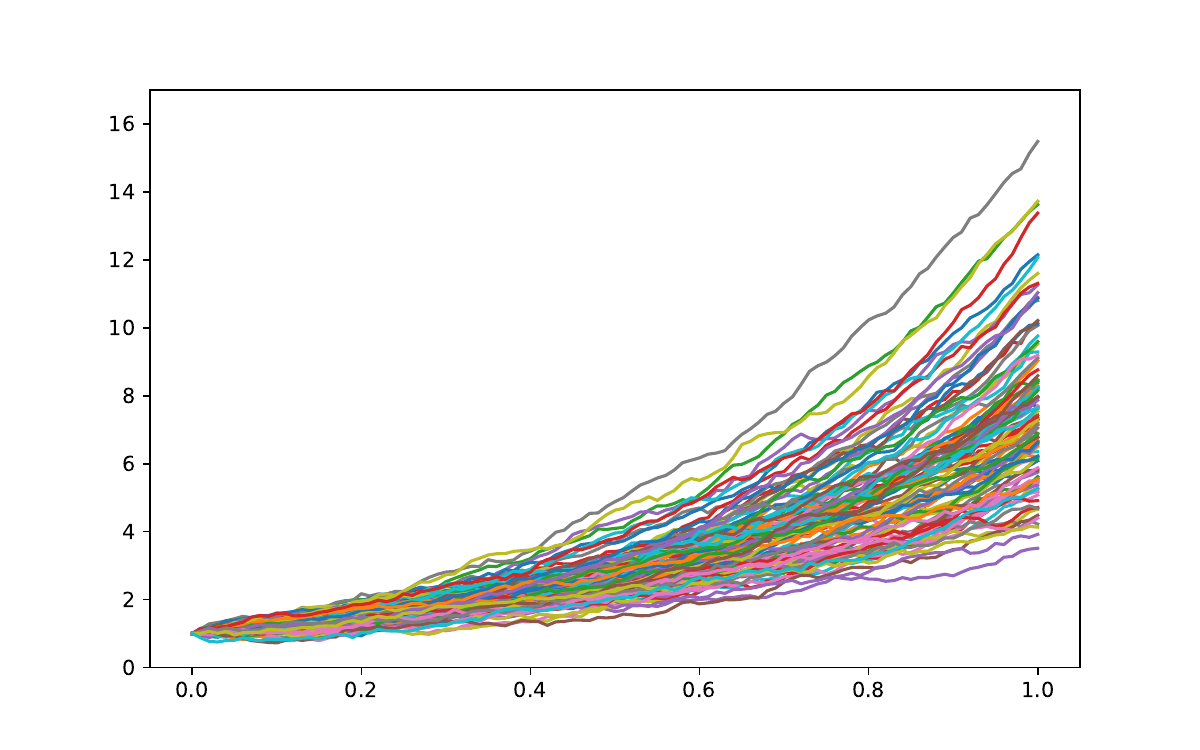}
    \par
    \includegraphics[width = 0.32\textwidth, trim=15mm 3mm 18mm 3mm,
    clip]{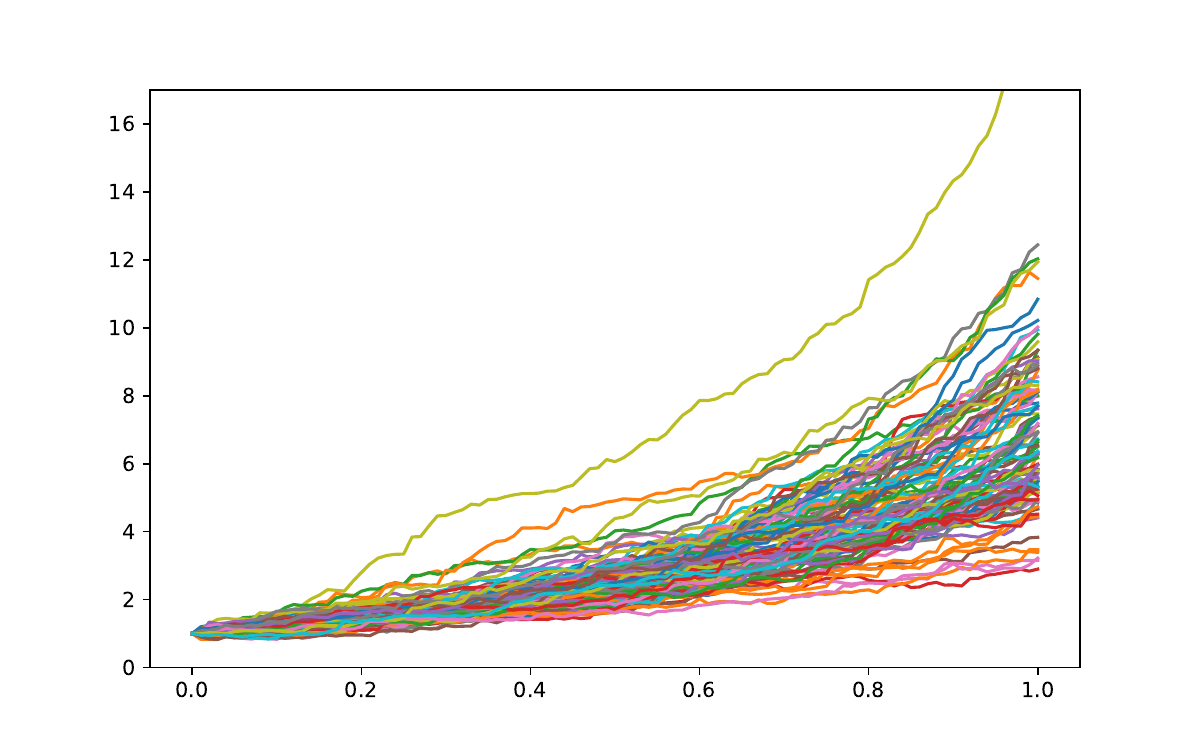}
    \includegraphics[width = 0.32\textwidth, trim=15mm 3mm 18mm 3mm,
    clip]{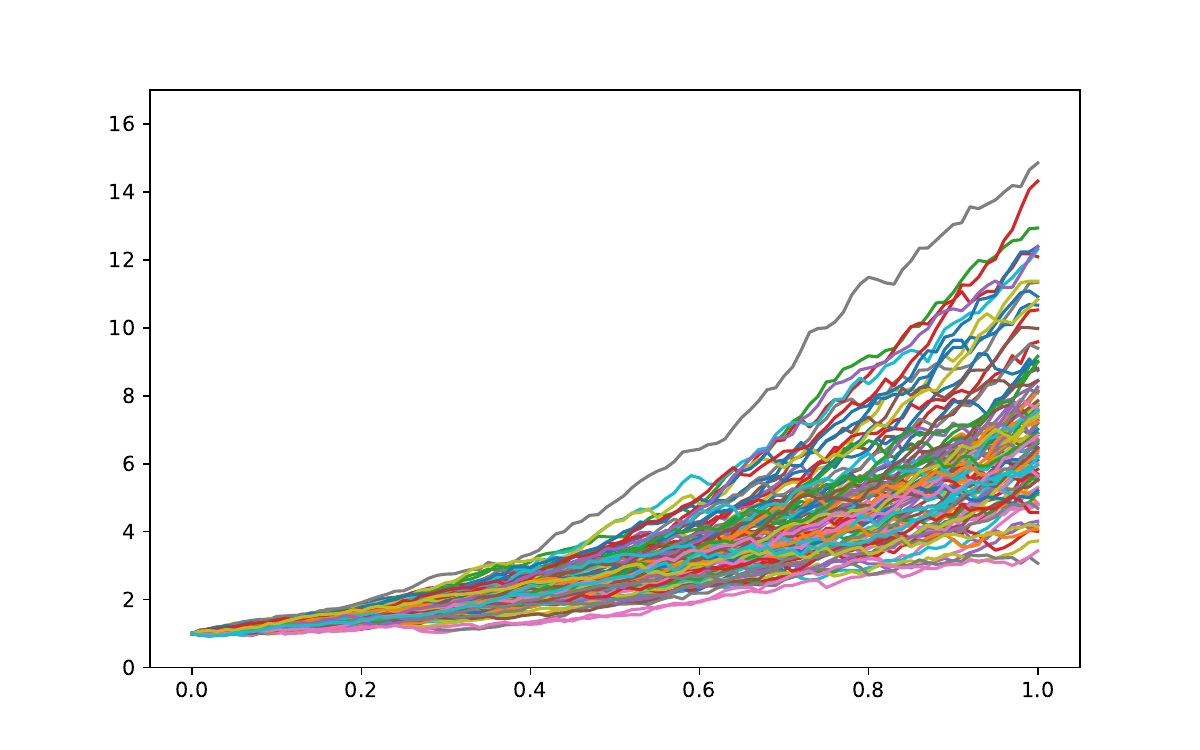}
    \includegraphics[width = 0.32\textwidth, trim=15mm 3mm 18mm 3mm,
    clip]{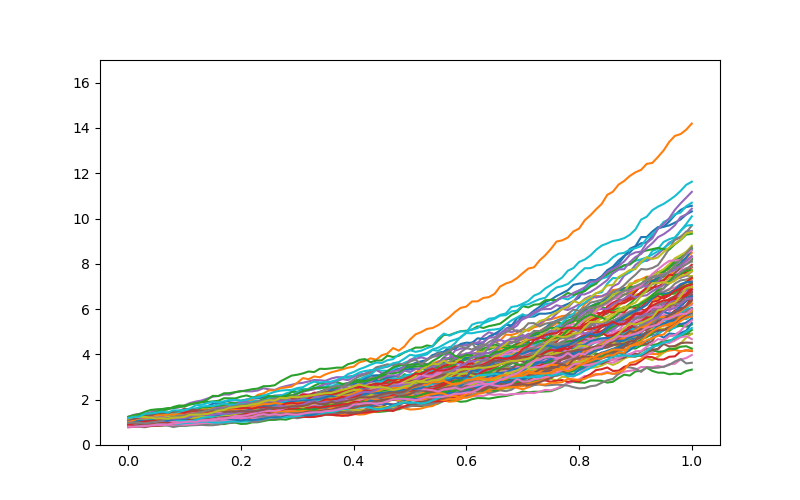}
    \caption{Trajectories generated by the following models trained in the 100\% observation regime. From left to right, top to bottom: JUMP, SDE,  MS ($\alpha=0.3$), TFM, SDE matching.}
\end{figure}

\paragraph{Additional Illustrations for Two-Dimensional Experiments}
Below, we provide trajectory plots corresponding to one experiment reported in Table \ref{tab:mmd_subsampling:stock2d}.

\begin{figure}[H]
    \centering
    \includegraphics[width=0.8\linewidth]{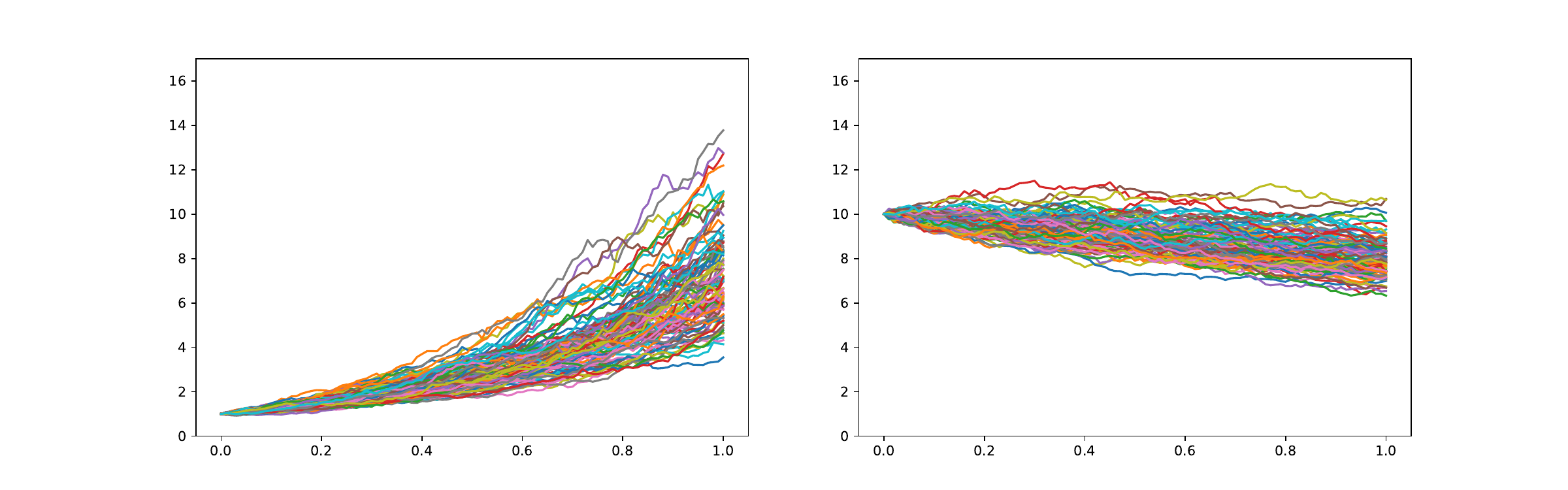}
    \caption{Ground truth trajectories of two-dimensional Black-Scholes dataset. Left: first coordinate. Right: second coordinate.}
\end{figure}

\begin{figure}[H]
    \centering
    
    \includegraphics[width = 0.49\textwidth, trim=15mm 5mm 15mm 5mm,
    clip]{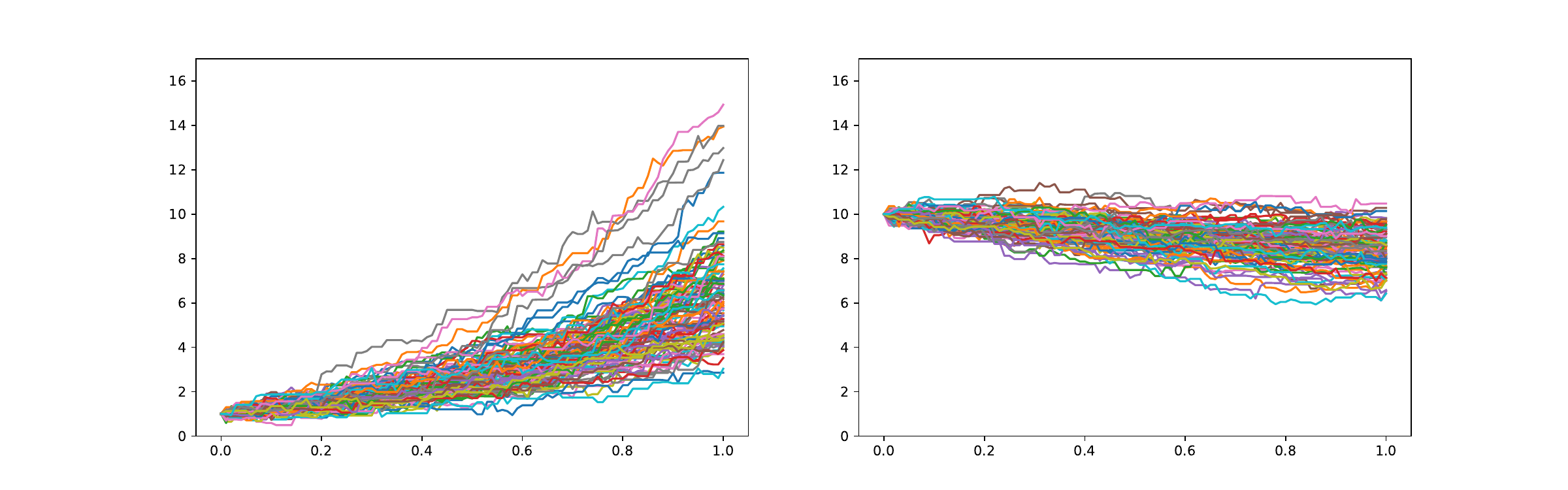}
    \includegraphics[width = 0.49\textwidth, trim=15mm 5mm 15mm 5mm,
    clip]{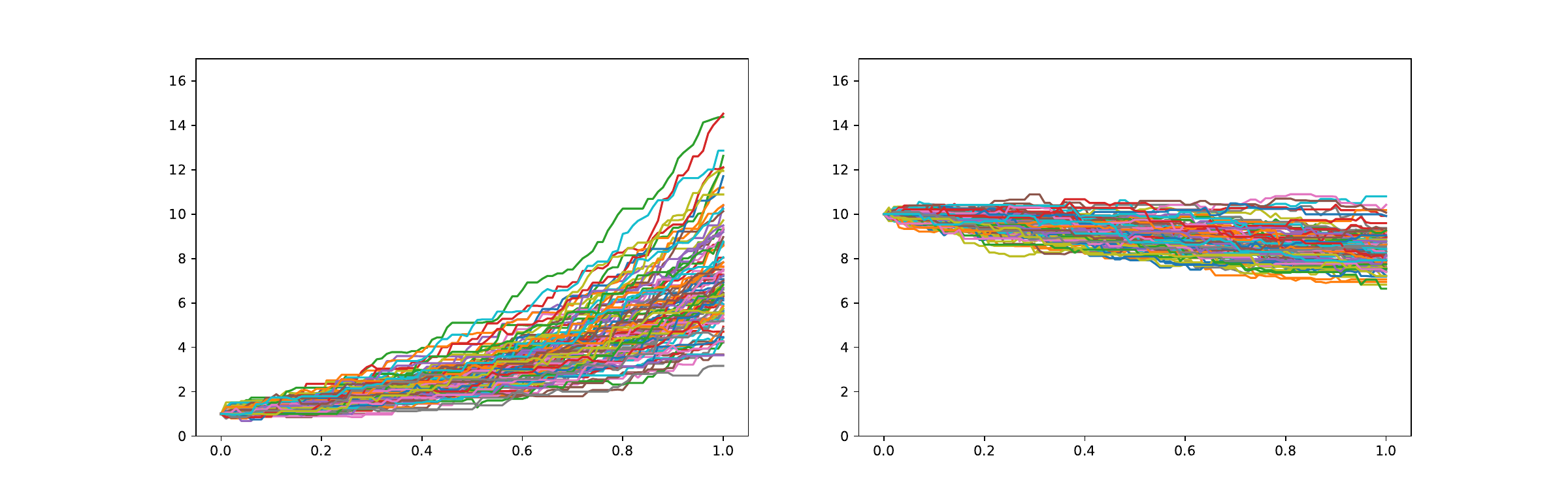}
    \par
    \includegraphics[width = 0.49\textwidth, trim=15mm 5mm 15mm 5mm,
    clip]{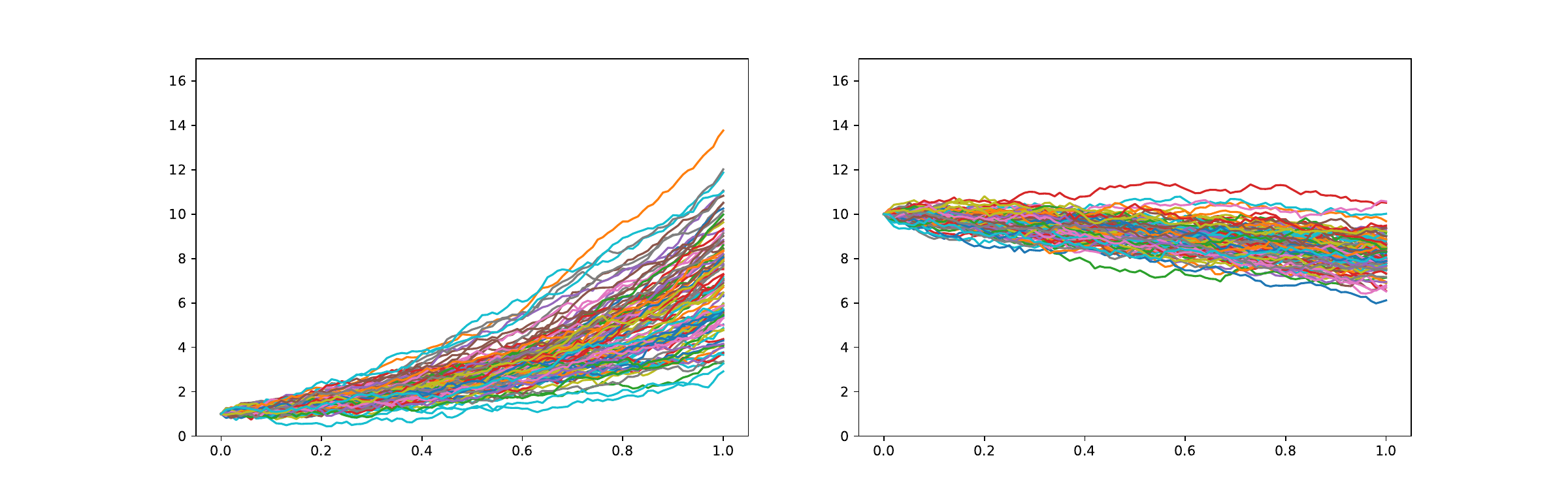}
    \includegraphics[width = 0.49\textwidth, trim=15mm 5mm 15mm 5mm,
    clip]{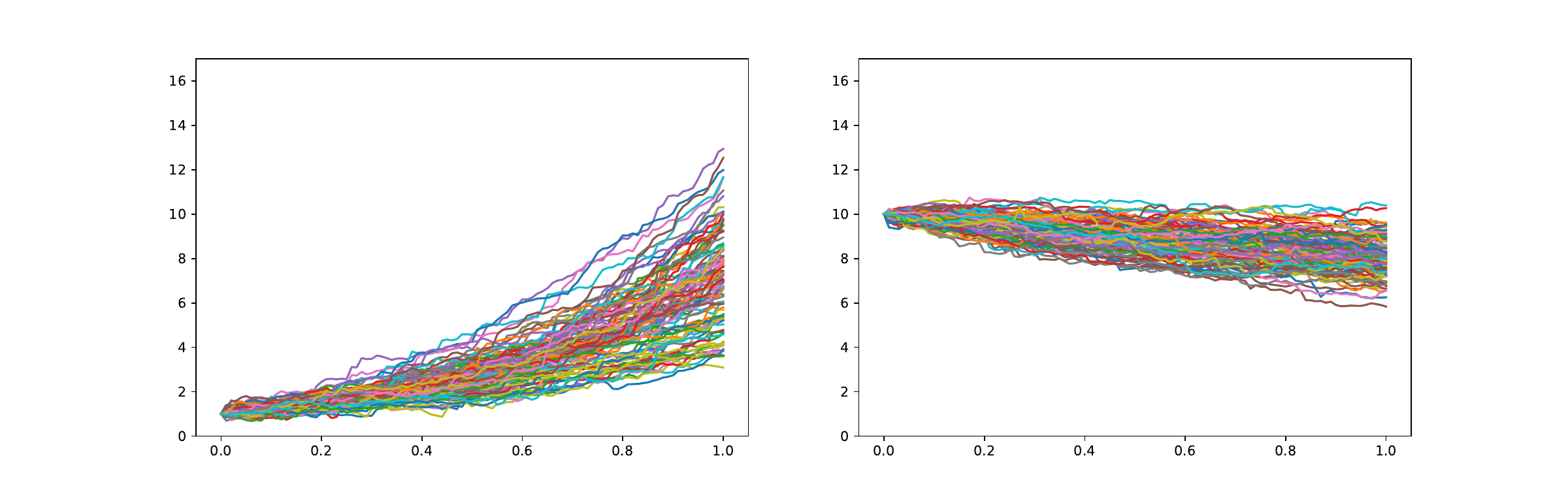}
    \par
    \includegraphics[width = 0.49\textwidth, trim=15mm 5mm 15mm 5mm,
    clip]{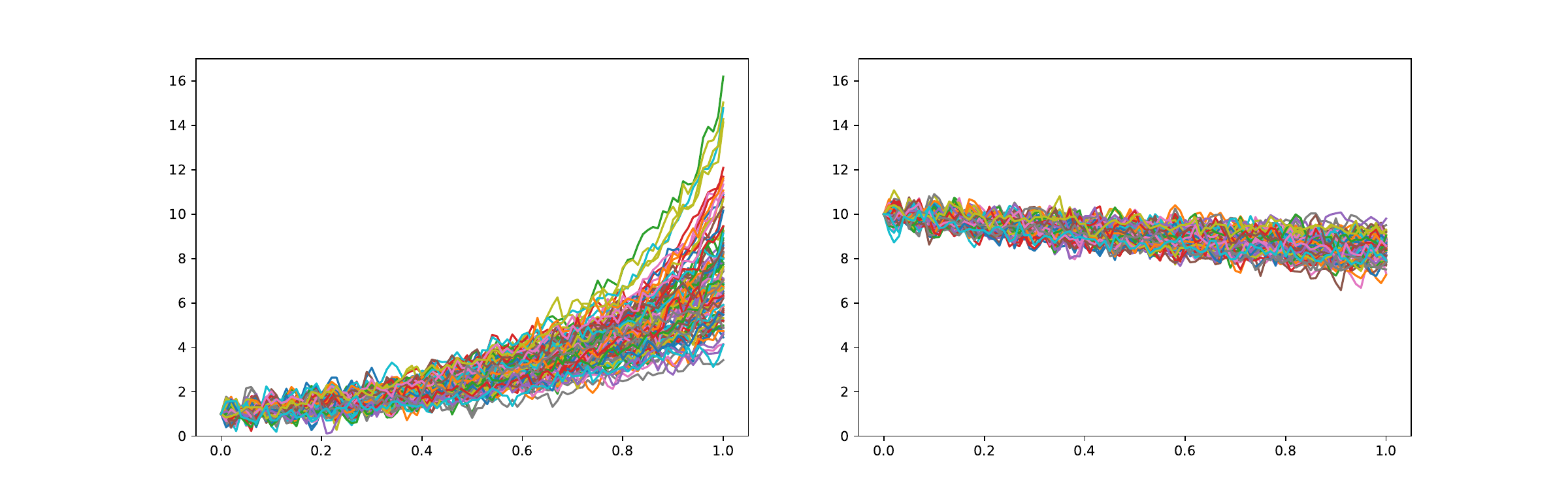}
    \phantom{\includegraphics[width=0.49\textwidth, trim=15mm 5mm 15mm 5mm,
    clip]{ima/images/stock2d/single/single_TfmModel_sub26_sigma10p0_seed0.pdf}}
    \caption{Illustration of '25\%'-column in Table \ref{tab:mmd_subsampling:stock2d}. Trajectories generated by the following models trained in the 25\% observation regime. From left to right, top to bottom: JUMP, JUMP (factorized components), SDE,  MS ($\alpha=0.4$), TFM.}
    \label{fig:Black-Scholes-2d-sub25}
\end{figure}

\section{Additional Details on the Experiment with Moving-Box Video Data}
\label{sec:supp:details_video}

The moving-box video dataset described in Section \ref{sec:video_data} of the article consists of binary videos of 30 frames on the time grid $\{0,1,\ldots,29\}$. We consider two different regimes, where either 20 or all of the 30 frames are observed. In the former case, we retain the first and last frame for each trajectory and select the remaining observed ones randomly.

We use $\rho^2 =0.01$, memory length $m=6$, and train only with full memory, i.e., during SGD on our loss functions $\mathcal{E}^{\textup{Diff}}$ and $\mathcal{E}^{\textup{Jump}}$, we only select times $t \geq t_{m-1}^i$ given a trajectory $(x^i_{t_0^i}, \ldots, x^i_{t_{n_i}^i})$.

\begin{figure}[htp]
    \centering
    \includegraphics[trim={38.3cm 22pt 26.7cm 4.4cm},clip,width=\linewidth]{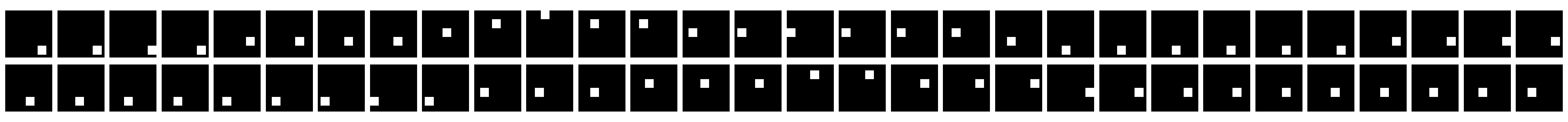}
    \includegraphics[trim={38.3cm 22pt 26.7cm 4.4cm},clip,width=\linewidth]{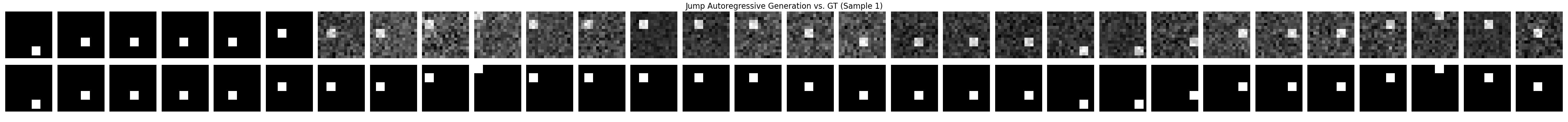}
    \includegraphics[trim={38.3cm 22pt 26.7cm 4.4cm},clip,width=\linewidth]{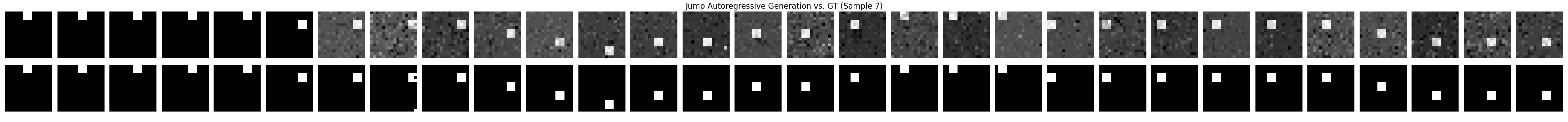}
    \caption{Learning of the dynamics of a moving box. Trajectories are given row-wise:  one ground truth trajectory (top row), one trajectory generated via a learned drift-diffusion generator (middle row) and via a learned jump generator (bottom row).}\label{moving:boxes}
\end{figure}

During autoregressive generation, the first $m$ frames of a fully observed ground truth trajectory are used as initial memory and the remaining frames are generated sequentially using interpolating bridges of length one. Recall that, during memory updates, the raw network output is binarized before it is appended to the memory: the nine largest pixel values are set to one, the remaining ones are set to zero. Figure \ref{moving:boxes} shows one ground truth and two generated trajectories.

While in the Black-Scholes experiments the bridge length during rollout is chosen to reflect the interval lengths observed during training, this is not sensible in this experiment due to the binary nature of the data. Observations at intermediate times between the endpoints of a bridge generally need not form a contiguous 3x3 box after binarization. Therefore, rollout with unit-length bridges is required, so that all observations occur at bridge endpoints. This is also the reason, why we do not consider a strongly subsampled regime: A strong mismatch between the interval lengths observed during training and the bridge lengths used for rollout is challenging for all compared methods.

We train each model in both subsampling regimes for ten different seeds and for $\eta^2 \in \{0.01, 0.1, 1\}$. As validation loss for checkpoint selection, we use the training objective on a set of held out validation trajectories. We define the box position given any selection of 9 pixels (not necessarily forming a contiguous 3x3 box) via the median $x$- and $y$-position. This way, the box position trajectory is defined for any rollout trajectory.

\begin{figure}[hbp]
    \centering
    \includegraphics[trim={38.3cm 0pt 26.7cm 23pt},clip,width=\linewidth]{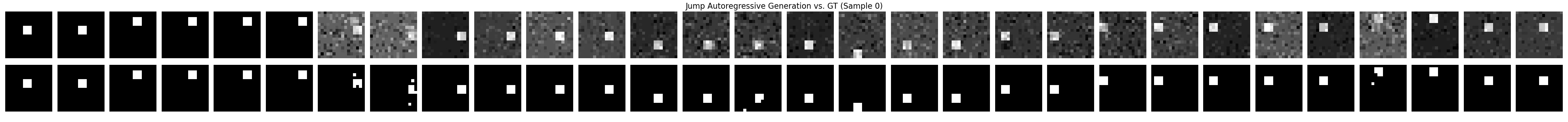}
    \caption{Partial rollout trajectory. Top: network output. Bottom: binarized output.}\label{fig:box_mov}
\end{figure}

In both regimes and for all models, we observe the minimal MMD for $\eta^2=1$. But, for every model trained with $\eta^2=1$, the majority of generated trajectories contains at least one frame, where the output does not show a 3x3 box, as illustrated in Figure \ref{fig:box_mov}. For models trained with $\eta^2=0.1$ or $\eta^2=0.01$, this almost never happened. For every model, at least 99\% of generated trajectories showed 3x3 boxes in all frames.

Therefore, we report in Table \ref{tab:mmd_subsampling:moving_box} the results for $\eta^2=0.1$, which are superior to the results for $\eta^2=0.01$ for all models. We do not report results for SDE matching  \cite{bartosh2025sde}, since this model is not well suited for capturing non-Markovian behavior (see the extended discussion on memoryless SDE-based models in \cite{zhang2024trajectory}).

\begin{table}
  \centering
  \begin{minipage}{\linewidth}
  \centering
  \scalebox{0.65}{
  \begin{tabular}{lll}
  \hline
  \textbf{Method} & \textbf{partially observed data} & \textbf{fully observed data} \\
  \hline
  Jump-based method
  & [0.0896 {\scriptsize$\pm$ 0.0347}]
  & [\textbf{0.0638} {\scriptsize$\pm$ 0.0186}] \\

  SDE-based method
  & [\textbf{0.0625} {\scriptsize$\pm$ 0.0213}]
  & [\textbf{0.0564} {\scriptsize$\pm$ 0.0120}] \\

  Jump + SDE (Markov superposition)
  & [\textbf{0.0515} {\scriptsize$\pm$ 0.0119}] {\scriptsize$(\alpha=0.9)$}
  & [\textbf{0.0662} {\scriptsize$\pm$ 0.0211}] {\scriptsize$(\alpha=0.8)$} \\

  TFM method \cite{zhang2024trajectory}
  & [0.1058 {\scriptsize$\pm$ 0.0180}]
  & [0.0916 {\scriptsize$\pm$ 0.0167}] \\
  \hline
  \end{tabular}
  }
  \vspace{1ex}
  \caption{Test MMD averaged over 10 seeds in the 20/30-frame and 30/30-frame observation regimes. All runs used $\eta^2=0.1$. For Markov superpositions, $\alpha$ denotes the SDE weight.}
  \label{tab:mmd_subsampling:moving_box}
  \end{minipage}
\end{table}

\end{document}